\documentclass[twocolumn,10pt]{autart}    
\usepackage{graphicx}          
\usepackage{amsmath,amssymb,mathrsfs}
\usepackage{amsmath}
\usepackage{enumerate}
\usepackage{graphicx}          
\usepackage{tikz}
\usepackage{cite}
\usepackage{natbib}

\usepackage{amsfonts}

\usepackage{amssymb}
\usepackage{subfigure}
\usepackage{hyperref}
\hypersetup{
	colorlinks = true,
	citecolor = cyan,
}
\usepackage{wrapfig}

\newtheorem{Theorem}{Theorem}[section]
\newtheorem{theorem}{Theorem}[section]

\newtheorem{lemma}[Theorem]{Lemma}
\newtheorem{proposition}[Theorem]{Proposition}

\begin{document}
\begin{frontmatter}

\title{Towards a Theoretical Foundation of PID Control for Uncertain Nonlinear Systems }
\date{}
\thanks[footnoteinfo]{
This work was supported by the National Natural Science Foundation of China under Grant No.11688101. (Corresponding author: Lei Guo)}

\author{Cheng Zhao}\dag\ead{zhaocheng@amss.ac.cn},
\author{Lei Guo}\dag\ead{lguo@amss.ac.cn},

\address{\dag Academy of Mathematics and Systems Science, Chinese Academy of Sciences, Beijing 100190, P.~R.~China}


\begin{abstract}
As is well-known, the classical PID control plays a dominating role in various control loops of industrial processes. However, a theory that can explain the rationale why the
linear PID can successfully deal with the ubiquitous uncertain nonlinear dynamical systems and a method that can provide explicit design formulae for the PID parameters are still lacking. This paper is a continuation of the authors recent endeavor towards establishing a theoretical foundation of PID for nonlinear uncertain systems.  In contrary to most of the existing literature on linear or affine nonlinear systems, we will consider a class of non-affine nonlinear uncertain systems, and will show that a three dimensional parameter set can be constructed explicitly, such that whenever the PID parameters are chosen from this set, the closed-loop systems will be globally stable and the regulation error will converge to zero exponentially fast, under some suitable conditions on the system uncertainties. Moreover, we will also consider the simpler PI and PD control, and  provide a necessary and sufficient condition for the choice of the PI parameters for a class of one dimensional non-affine uncertain systems, by applying the Markus-Yamabe theorem in differential equations. These theoretical results show explicitly that  the ubiquitous PID control does indeed have strong robustness with respect to both the system nonlinear uncertainties and the selection of the controller parameters.
\end{abstract}
\end{frontmatter}

\section{Introduction}
Feedback has had a revolutionary influence in practically all areas it is used, and the classical proportional-integral-derivative (PID) control is perhaps the most basic form of feedback (\cite{As}). In fact, the PID controller is by far the most basic and widely used control algorithm in engineering systems. For example, more than 95\% of the control loops in process control are of PID type, where most loops
are actually PI control (\cite{Astrom19955}). It is also reported that the PID control has a much higher influence than the advanced control technologies, and that we still have nothing that compares with PID by far (\cite{samad}).

Therefore, a natural question is: why the linear PID performs so successfully in controlling the ubiquitous nonlinear systems with various uncertainties. There are some partially known reasons: It has a simple structure whose design is model-free; it has the ability to reduce the influence of uncertainties through proportional action, to eliminate steady state offsets through integral action and to anticipate the tendency of the output through the derivative action; the  Newton's second law is fundamental in modeling mechanical systems, which is just suitable for the PID control. However, these reasons appear to be superficial, and a fundamental theory that can guarantee the global stability and strong robustness of the classical PID when applied to  uncertain nonlinear dynamical systems together with an explicit design method for the PID parameters are still lacking.

Traditionally, the PID parameters are chosen based
on experiments or experiences or by both in practice, including the well-known Ziegler-Nichols rules (\cite{ziegler1942}). However, it has been reported that most of the practical PID loops are poorly tuned, and there is strong evidence that PID controllers remain poorly understood (\cite{o2006pi}). This may partly due to the fact that the PID controller has not attracted enough attention from the research community (\cite{Astrom19955}). In fact, most related theoretical studies of PID are conducted for systems with special structures, e.g., for  linear systems(e.g., \cite{As,Astrom19955,Ho2003,Keel2008,Silva}), for affine nonlinear systems \cite{Khalil,Killingsworth2006}. There are also abundant works focused on nonlinear robot systems, see e.g., \cite{Takegaki,Romero,Alvarez,Besancon,Su,Borja}, etc. A notable work \cite{Takegaki} in robotics  proved that a simple PD law provides a global solution to the point-to-point positioning task for fully actuated robot manipulators. The global asymptotic  regulation of uncertain robotic  systems via the PD control plus a class of nonlinear integral action has been established by \cite{kelly}. It is worth noting that for robot systems, the total energy stored in the systems decreases with time and the dissipativity properties  can be used to design stable and robust feedback controllers. Moreover, the stability of a dissipative system with a PID controller can be established using dissipative arguments which hinges  on the specific structure of the robot systems (\cite{Brogliato}).

However, for more general nonlinear uncertain systems without dissipative properties, it is more challenging to analyze the asymptotic behavior of the control system since there is no concept like the total energy  which can be used as a natural Lyapunov function candidate. The longstanding gap between the PID theory and its widespread practice calls for a theory for general nonlinear uncertain control systems (\cite{lguo}).
These motivate our recent studies on the theory and design of the PID control for nonlinear uncertain systems(see e.g., \cite{Zhao2017,Zhang2019,zhao2020}). For example, we have provided some global convergence results for a class of second order uncertain nonlinear systems in our previous works \cite{Zhao2017,Zhang2019} with no uncertainty in the control channel. For affine-nonlinear uncertain systems with a general relative degree,  we have investigated the capability of the extended PID control in \cite{zhao2020}, where both the input and output are assumed to be one dimensional. Undoubtedly, most mechanical systems such as aircrafts, manipulators or walking robotics have more than one degree of freedom in the control inputs. Moreover, there are many practical systems are in the non-affine from, such as flight control systems \cite{bos} , and chemical reactions, etc, or systems with nonsymmetric control gain matrix, see e.g., \cite{huang2001}.
These facts inspired us to consider general MIMO non-affine uncertain systems based on our previous works (\cite{Zhao2017,Zhang2019}). It is worth mentioning that both the controller design and analysis of MIMO uncertain systems considered in this paper are more challenging than those studied previously, because of the uncertainties caused by the strong coupling of the MIMO non-affine nonlinear control systems.

The main contributions are as follows.  First,  we will show that for  a basic class of multi-input multi-output(MIMO) second order non-affine uncertain systems,  a 3-dimensional set can be constructed from which the three PID parameters can be chosen arbitrarily to globally stabilize the closed-loop control system and to make the regulation error converge to zero exponentially fast, as long as some certain bounds of the partial derivatives of the uncertain function are known a prior. Next, for the case where the setpoint is an equilibrium of the uncontrolled system or the case where the control system is of first order, we will show that the PD (or PI) controller can globally regulate the considered systems. Finally, we will use the  Markus-Yamabe theorem to establish a necessary and sufficient condition for the choices of PI parameters for a class of one dimensional nonlinear uncertain systems.  These results will demonstrate explicitly that the PID control does indeed have large scale robustness with respect to both the uncertain system structure and the selection of the controller parameters, where the controller gains are not necessarily to be high.

The rest of the paper is organized as follows. The problem formulation will be described  in  Section 2. Section $3$ will present our  main results, with proofs put in Section $4$. Section $5$ will conclude the paper with some remarks.

\section{Problem Formulation}
\subsection{Notations}
Denote $\mathbb{R}^n$ as the $n$-dimensional Euclidean space, $~\mathbb{R}^{m\times n}$ as the space of $m\times n$ real matrices, $\|x\|$  as the Euclidean norm of a vector $x$, and $x^{\mathsf{T}}$ as the transpose of a vector or matrix $x$.
 The  norm of a matrix $P \in \mathbb{R}^{m\times n}$  is defined by $\|P\|=\sup_{x\in\mathbb{R}^n, \|x\|=1} \|Px\|$. For a square matrix $P\in\mathbb{R}^{n\times n}$,  denote $\text{Sym} [P]\overset{\triangle}{=}(P+P^{\mathsf{T}})/2$ as the symmetrization matrix of $P$. For a symmetric matrix $S\in \mathbb{R}^{n\times n}$, we denote $\lambda_{\min}(S)$ and $\lambda_{\max}(S)$ as the smallest and the largest eigenvalues of $S$, respectively. For two symmetric matrices $S_1$ and $S_2$ in $\mathbb{R}^{n\times n}$, the notation $S_1>S_2$ or $S_2<S_1$ implies that  $S_1-S_2$ is a positive definite matrix; $S_1\geq S_2$ or $S_2\le S_1$ implies  that $S_1-S_2$ is a positive semi-definite matrix.

Let $\mathbb{R}_+^2$ denote the product $(0,\infty)\times(0,\infty)$,  and similarly  $\mathbb{R}_+^3=(0,\infty)\times(0,\infty)\times(0,\infty)$. Let $C^{1}(\mathbb{R}^{n},\mathbb{R}^{m})$ be the space of continuously differentiable functions from $\mathbb{R}^{n}$ to $\mathbb{R}^m$, denoted as   $C^{1}(\mathbb{R}^{n})$ for simplicity when $m=1$.
For a function $\Phi=(\phi_1,\cdots,\phi_m)^{\mathsf{T}} \in C^{1}(\mathbb{R}^{n},\mathbb{R}^{m})$, let
$$
\frac{\partial\Phi}{\partial x}(x)=\begin{bmatrix}\begin{array}{cccc}
\frac{\partial\phi_{1}}{\partial x_{1}}(x)~ & \frac{\partial\phi_{1}}{\partial x_{2}}(x)& \cdots & ~ \frac{\partial\phi_{1}}{\partial x_{n}}(x)\\
		\vdots & \vdots &\ddots &  \vdots\\
\frac{\partial\phi_{m}}{\partial x_{1}}(x)~ &\frac{\partial\phi_{m}}{\partial x_{2}}(x) &\dots & ~ \frac{\partial\phi_{m}}{\partial x_{n}}(x)
\end{array}\end{bmatrix}$$ denote the Jacobian of $\Phi$ at the point $x$.

A function $\Phi \in C^1(\mathbb{R}^n,\mathbb{R}^n)$ is called a global diffeomorphism on $\mathbb{R}^n$ if it is both \emph{injective} and  \emph{surjective}, and  its inverse satisfies  $\Phi^{-1}\in C^1(\mathbb{R}^n,\mathbb{R}^n)$.

\subsection{The control system}
Consider the following class of MIMO non-affine uncertain nonlinear systems:
\begin{equation}
\begin{cases}
\dot{x}_{1}&=x_{2}\\
\dot{x}_{2}&=f(x_{1},x_{2},u),~~x_1\in\mathbb{R}^n,~u\in\mathbb{R}^n,\\
\end{cases}\label{equation}
\end{equation}
where $(x_1,x_2)$ is the system state that can be measured, $u$ is the control input and $f\in C^1(\mathbb{R}^{3n}, \mathbb{R}^{n})$ is an \emph{uncertain nonlinear} function.

We remark that many practical dynamical systems can be described by the basic model (1) via the Newton's second law. For example,  it can  be used to describe the motion of a controlled moving body in $\mathbb{R}^n$, where $x_1$ and $x_2$ represent the position and velocity of the moving body respectively and where the nonlinear function $f(x_{1},x_{2},u)$ represents the total  external forces acting on the moving body. It can also be used to model a general mechanical system with $n$ degrees of freedom, where $x_1$ and $x_2$ represents the generalized coordinates and the generalized velocity of the system respectively.

Our control objective is  to globally stabilize the system (\ref{equation}) and  to make the controlled variable $x_1(t)$ converge to a desired reference value $y^{*}\in\mathbb{R}^n$ exponentially for all initial states $(x_1(0),x_2(0))\in\mathbb{R}^{2n}$, under the condition that  the nonlinear function $f(\cdot)$ contains uncertainty.
\section{The Main Results}
\subsection{PID control}
In this section, we will investigate the capability of the
classical PID controller:
\begin{align}
u(t)=&k_{p} e(t)+k_{i} \int_{0}^{t} e(s)ds+k_{d} \dot{e}(t),\\
e(t)=&y^{*}-x_1(t),\nonumber
\end{align}
 and provide  a concrete design method for the three PID parameters $(k_p,k_i,k_d)$ to achieve the desired objectives.

Note that $f$ is  uncertain, to establish a rigorous mathematical theory, we need to introduce a measure on the size of the class of uncertain functions. It is conceivable that  the sensitivities of the uncertain function with respect to their variables, reflected by the upper bounds of the partial derivatives  $\frac{\partial f}{\partial x_i}$($i=1,2$), will play an essential role in the control of uncertain systems,  see (\cite{xie2000,Zhao2017}). In general, such upper bounds may be unbounded functions of the state variables, in which cases one can only get semi-global control results(cf. e.g., \cite{Khalil,zhao2020}) if the open-loop system is unstable and no additional structural properties such as dissipativity are available, because of the linear structure of the classical PID. Moreover, in order for the input signal to have the ability to influence the systems to be controlled, the control gain matrix $\frac{\partial f}{\partial{u}}$ should not vanish. These natural intuitions and facts inspired us to introduce the  following function space.
\begin{defn}[Function space]
For given positive constants $L_1$, $L_2$ and $\underline{b}$, define $\mathcal{F}(=\mathcal{F}_{L_1,L_2,\underline{b}})$ as follows:
 \begin{align*}
 \mathcal{F}\overset{\triangle}{=}\!
 \bigg \{f\!\in \!C^1(\mathbb{R}^{3n}, \mathbb{R}^{n})~\!\bigg|~\!
 \Big\|\frac{\partial f}{\partial{x_1}}&\Big\|\le L_1,~\Big\|\frac{\partial f}{\partial{x_2}}\Big\|\le L_2,\\
 {\rm Sym} \Big[~\!\frac{\partial f}{\partial{u}}~\!&\Big]\geq \underline{b}I_n,~ \forall x_1,x_2,u\in\mathbb{R}^n\bigg\},
 \end{align*}
where $\frac{\partial f}{\partial{x_1}}$, $\frac{\partial f}{\partial{x_2}}$ and $\frac{\partial f}{\partial{u}}$ are the $n\times n$ Jacobian matrices of $f$ with respect to $x_1,x_2$ and $u$, respectively, $I_n$ is the $n\times n$ identity matrix and
${\rm Sym}\big[\frac{\partial f}{\partial{u}}\big]\overset{\triangle}{=}\big[\frac{\partial f}{\partial{u}}+\big(\frac{\partial f}{\partial{u}}\big)^{\mathsf{T}}~\!\big]\big/2.$
\end{defn}
We now give some explanations for the uncertain function space $\mathcal{F}$.  First, it is worth mentioning that for nonlinear functions $f\in\mathcal{F}$, the uncontrolled system (\ref{equation}) may be unstable. Moreover,  no specific structural information is assumed in this paper except for some prior information on the  bounds of the partial derivatives. Since the classical PID is a linear feedback, for systems (1)  without dissipative properties, the boundedness of the partial derivatives $\frac{\partial f}{\partial x_i}$($i=1,2$)(or the linear growth condition) appears to be necessary in general  for global results (\cite{zhao2016}). Of course, the linear growth condition can be considerably relaxed if we are only interested in obtaining semi-global results (\cite{zhao2020}), and one may even get global results if specific structural information such as that in robot manipulators is available (\cite{kelly}).
Furthermore, we remark that the constants $L_1$ and $L_2$ correspond to the upper bounds of the ``anti-stiffness" and the ``anti-damping"  coefficients of the nonlinear systems (1) respectively (\cite{krstic2017}), and that the three positive constants $L_1$, $L_2$ and $\underline b$
 can be used to describe the system uncertainty quantitatively, since the ``size'' of the uncertain function space $\mathcal{F}_{L_1,L_2,\underline{b}}$  will increase when either $L_1$ and $L_2$ increases or $\underline b$ decreases. Finally, we mention that the meanings of the symbols $\frac{\partial f}{\partial x_i}$ and $\frac{\partial f}{\partial u}$ remain the same as in Definition 1 throughout the paper.

Based on the triple $(L_1,L_2,\underline b)$, we next introduce a three dimensional unbounded parameter set.
\begin{defn}[Parameter set] Define $\Omega_{pid}$  as follows:
\begin{align}\label{pid}
\Omega_{pid}=\left\{(k_p,k_i,k_d)\in\mathbb{R}_+^3 \left|
\begin{array}{c}
		~k_p^2>2k_ik_d+\bar k\\
		~k_d^2~\!>k_p/\underline b~\!+\bar k
	\end{array}\right.
\right\},
\end{align}
where $\bar k\overset{\triangle}{=}(L_1+L_2)(k_p+k_d)/\underline b$.
\end{defn}

\textbf{Remark 3.1}~~We remark that for any given positive constants $L_1,L_2$ and $\underline b$, the parameter set $\Omega_{pid}$  defined by (\ref{pid}) is an open and unbounded subset in $\mathbb{R}^3$.  In fact,  for any given $k_i>0$, it can be  verified that $(k_p,k_i,k_d)\in\Omega_{pid}$, provided that
$$ k_p=k_d\geq 2k_i+[2(L_1+L_2)+1]/\underline b.$$
Hence, the integral parameter $k_{i}$ can be chosen arbitrarily small, and the PD gains $(k_p,k_d)$ need not to be sufficiently large.  Moreover, the parameter set $\Omega_{pid}$  is actually a semi-cone in the sense that:
$$(k_p,k_i,k_d)\in\Omega_{pid}\Rightarrow\alpha(k_p,k_i,k_d)\in\Omega_{pid},~~~\forall \alpha\geq 1.$$ The proof of Remark 3.1 is given in \S \ref{re1}.

Now, we are in position to present the main results.
\begin{theorem}
Consider the PID controlled system (1)-(2). Then for any $(k_p,k_i,k_d)\in \Omega_{pid}$, there exist two constants $M>0$ and $\lambda>0$(depend on $(k_p,k_i,k_d,L_1,L_2,\underline b)$ only), such that for any $f\in \mathcal{F}_{L_1,L_2,\underline{b}}$,~any setpoint $y^*\in \mathbb{R}^n$ and any initial states $(x_1(0),x_2(0))\in \mathbb{R}^{2n}$, the solution of the closed-loop system will satisfy
\begin{align*}\|e(t)\|+\|\dot e(t)\|\le
Me^{-\lambda t} \left(\|e(0)\|+\|\dot e(0)\|+\|u^*\|\right)\end{align*}
for all $t\geq 0$, where $e(t)=y^*-x_1(t)$ and $u^*\in\mathbb{R}^n$ is the unique solution of the algebraic equation $f(y^*,0,u)=0$.
\end{theorem}
The proof of Theorem 3.1 is  given in \S \ref{pf}.

\textbf{Remark 3.2} ~From Theorem 3.1, we know that the classical PID (2) can globally stabilize and regulate the uncertain system (1), provided that the upper bounds of  $\frac{\partial f}{\partial x_i}$($i=1,2$) together with a lower bound to $\frac{\partial f}{\partial u}$ are known \emph{a prior}. It is worth mentioning that both the upper bound and the symmetry of the control gain matrix $\frac{\partial f}{\partial u}$ are not required,  and $x_1$ is allowed to be high-dimensional. Thus, Theorem 3.1 has weakened the assumptions on the system nonlinear functions considered in \cite{Zhang2019} and \cite{Zhao2018} significantly, where \cite{Zhang2019} only considers affine-nonlinear systems with no uncertainty in the control channel, and  \cite{Zhao2018} only considers the one-dimensional case. In addition, notice that the unbounded set $\Omega_{pid}$ only depends on $L_1,L_2$ and $\underline{b}$, we can see that  the selection of the PID parameters relies neither on the precise structure of the function $f(\cdot)$, nor on the initial states and  the setpoint $y^*$. Hence, the PID controller has strong robustness with respect to both the system uncertainties and  to the selection of the three controller parameters. This partially explains the rationale behind the widespread successful applications of the PID control.

\textbf{Remark 3.3}
~Theorem 3.1 is also true for bounded time-varying reference signals $y^*(t)$, save that the output tracking error may not converge to zero in general,  which can be made arbitrarily small by either assuming slowly varying reference signals,  or choosing large PID controller gains, see e.g.,  \cite{jiang2015} for related results with extended state observer-based controller. For the case where only the state variable $x_1$ is measured, a high-gain differential observer may be used to estimate the derivative signal to get the similar results, see \cite{guo2021}  in the affine case.

\textbf{Remark 3.4} ~~ It is possible to further specify the choice of the three PID parameters $(k_p,k_i,k_d)$ from the set $\Omega_{pid}$. One way is to rewrite the PID control as an adaptive pole-placement control  where the uncertain nonlinear function is adaptively ``cancelled" by an online estimator formed by a suitable linear combination of the three terms in PID, see \cite{Zhong} for a related discussion based on an earlier version of the current paper.

\textbf{Remark 3.5} ~~In Theorem 3.1, the partial derivatives $\frac{\partial f}{\partial x_i}, i=1,2$ of the nonlinear function $f(x_1,x_2,u)$ have been assumed to be bounded for establishing global results. Although it seems to be somewhat restrictive mathematically and may not be satisfied for some specific cases including robot manipulators(see e.g., \cite{Besancon,kelly}), it is a quite natural assumption in general. Firstly, the boundedness of the partial derivatives is equivalent to the standard Lipschitz condition used in ordinary differential equation to ensure the global existence and uniqueness of solutions. Secondly, for sampled-data control of nonparametric nonlinear uncertain systems, the boundedness of the derivatives of the uncertain functions appears to be necessary for global stabilization(see \cite{xue2002,xie2000,lguo}).  Thirdly, this condition can be verified for many practical systems since the state variables are usually kept in a bounded set due to physical constraints, and  even when there were no physical constrains on the state variables, it may also be verified for many practical systems, e.g., damped vibration,  simple pendulum with damping, etc.

\subsection{PD control}
It is known that when the setpoint is an \emph{equilibrium} of the uncontrolled system, then the integral term  is not necessary for regulation for a class of affine nonlinear system with $n=1$, see Theorem 2 in \cite{Zhao2017}.

For non-affine uncertain system (\ref{equation}), let us consider the case  where the setpoint $y^*\in\mathbb{R}^n$ is an \emph{equilibrium} of the uncontrolled system, i.e., $f(y^*,0,0)=0$, and introduce the following function space:
\begin{align*}
\mathcal{F}_{L_1,L_2,\underline{b},y^*}\overset{\triangle}{=}
\big\{f\in\mathcal{F}_{L_1,L_2,\underline{b}}~ \big|~f(y^*,0,0)=0\big\}.
\end{align*}

Define the following two dimensional PD parameter set:
 \begin{align}\label{pd1}
 \Omega_{pd}=\Big\{ (k_p,k_d) \in \mathbb{R}_+^2 ~\Big| ~\!~k_p^2>\bar k, ~k_d^2>k_p/\underline b+\bar k~\!\Big\},
\end{align}
where $\bar k\overset{\triangle}{=}(L_1+L_2)(k_p+k_d)/\underline b$.
\begin{theorem}
Consider the  non-affine uncertain system (1), where $u(t)$ is the following PD control:
\begin{equation}
u(t)=k_{p} e(t)+k_{d} \dot{e}(t), ~~~e(t)=y^{*}-x_1(t).
\end{equation} Then for any $(k_p,k_d)\in \Omega_{pd}$, there exist constants $M>0$ and $\lambda>0$(depend on $(k_p,k_d,L_1,L_2,\underline b)$ only), such that  the closed-loop system will satisfy
\begin{equation*}\|e(t)\|+\|\dot e(t)\| \le  Me^{-\lambda t} \big(\|e(0)\|+\|\dot e(0)\|\big),\end{equation*}
for all $t\geq 0$, for all $f\in \mathcal{F}_{L_1,L_2,\underline{b},y^*}$ and for all initial states $(x_1(0),x_2(0))\in \mathbb{R}^{2n}$.
\end{theorem}

\textbf{Remark 3.5}~~It can be  verified that $(k_p,k_d)\in\Omega_{pd}$ provided that  $k_p=k_d=k>[2(L_1+L_2)+1]/\underline b,$ since
\begin{align*}
k_p^2-\bar k=&k_d^2-\bar k >k_d^2-k_p/\underline b-\bar k\\
=&k^2-k/\underline b-2k(L_1+L_2)/\underline b\\
=&k(k-[2(L_1+L_2)+1]/\underline b)>0.
\end{align*}
Thus, $\Omega_{pd}$ is an open and unbounded subset in $\mathbb{R}^2$ and the PD gains $(k_p,k_d)$ need not to be sufficiently large to satisfy $(k_p,k_d)\in\Omega_{pd}$. Moreover, we remark that in Theorem 3.2, if the dimension of $x_1$ is one, i.e., $n=1$,  then a necessary and sufficient condition for the choice of PD parameters can be found, see \cite{Zhao2018}.

\subsection{PI control}
In this subsection, we consider the following class of first order non-affine nonlinear uncertain systems:
\begin{equation}
\dot x=f(x,u),~~x\in\mathbb{R}^n,~u\in\mathbb{R}^n,\\
\label{equation2}
\end{equation}
where $u$ is the control input and $f(x,u):\mathbb{R}^{2n}\to \mathbb{R}^{n}$ is an \emph{uncertain} nonlinear function.  It is worth noting that (\ref{equation2}) can be used to model many industrial processes, for examples, the level controls in single tanks, stirred tank reactors with perfect mixing, etc, see \cite{Astrom19955}.

Next, we will show that the following PI controller:
\begin{equation}\label{pi}u(t)=k_{p} e(t)+k_{i} \int_{0}^{t} e(s)ds,~~ ~e(t)=y^{*}-x(t)\end{equation}
can globally stabilize and regulate the system (\ref{equation2}).

Consider the following class of uncertain functions:
\begin{align*}
\mathcal{F}_{L,\underline{b}}\overset{\triangle}{=}
\bigg \{f\in C^{1} \left|~
 \big\|\frac{\partial{f}}{\partial{x}}\big\|\le L,~
\text{Sym}\Big[\frac{\partial f}{\partial u}\Big]\geq \underline{b}I_n,
\forall x,u\right.\bigg\},
\end{align*}where $L>0$ and $\underline{b}>0$ are two constants.
Define the following two dimensional PI parameter set:
\begin{align*}
\Omega_{pi}=\Big\{( k_p,k_i)\in \mathbb{R}_+^2~\Big|~ k_p^2\underline b>k_pL+k_i+L^2/(4\underline b)\Big\}.
\end{align*}
\textbf{Remark 3.6}~~It can be  verified that $(k_p,k_i)\in\Omega_{pi}$ provided that $k_i>0$ and $k_p\geq 2L/\underline b+k_i/L$, since
\begin{align*}
k_p^2\underline b-k_pL-k_i=&k_p(\underline b k_p-L)-k_i\geq k_pL-k_i
\\\geq &2L^2/\underline b+k_i-k_i >L^2/(4\underline b).
\end{align*}
Thus, $\Omega_{pi}$ is an open and unbounded subset in $\mathbb{R}^2$.

\begin{theorem}
Consider the uncertain nonlinear system (\ref{equation2}), where $u(t)$ is the  PI controller defined by (\ref{pi}).
Then whenever $(k_p,k_i)\in \Omega_{pi}$, the closed-loop system will satisfy
$\lim_{t\to\infty} x(t)=y^*$,
for any $f\in \mathcal{F}_{L,\underline{b}}$, any setpoint $y^*\in \mathbb{R}^n$ and any initial state $x(0)\in \mathbb{R}^{n}$.  In addition, if $\dim x=\dim u=1$, then the PI parameters can be chosen from the following larger and necessary parameter set:
\begin{align*}
\Omega'_{pi}=\big\{( k_p,k_i)\in \mathbb{R}^2~\big | ~k_p\underline b>L,~~k_i>0\big\}.
\end{align*}
\end{theorem}

\section{Proof of the Main Results}\label{pf}
We first list some auxiliary results  that will be used in the proofs of the main results.

The following result is known as Hadamard's Theorem, which asserts that smooth($C^1$) and proper maps with non-singular Jacobian are
global diffeomorphisms. To be specific,

\textbf{Theorem A1} (\cite{gordon1972,Ruzhansky2015}) Let $\Phi \in C^1(\mathbb{R}^n,\mathbb{R}^n)$, then $\Phi$ is a global diffeomorphism on $\mathbb{R}^n$ if and only if the Jacobian matrix $\frac{\partial \Phi}{\partial x}$ is nonsingular for all $x\in \mathbb{R}^n$ and $\lim_{\|x\|\to \infty}\|\Phi(x)\|=\infty$.

\textbf{Theorem A2} (\cite{fessler1995proof}) Consider the ordinary  differential equation $\dot{x}=f(x),$ where  $f\in C^1(\mathbb{R}^2,\mathbb{R}^2)$ and $f(0)=0$. If for any $x\in \mathbb{R}^2$, the eigenvalues of the  Jacobian $Df(x)$  have negative real parts, then the zero solution of the differential equation is globally asymptotically  stable.

\begin{proposition}
Let $\Phi \in C^1(\mathbb{R}^n,\mathbb{R}^n)$. Suppose that  \begin{align*}{\rm Sym}\Big[\frac{\partial \Phi}{\partial x}\Big]\geq R,~~\forall x\in\mathbb{R}^n,\end{align*}
where $R\in\mathbb{R}^{n\times n}$ is a constant positive definite matrix,  then  $\Phi$ is a global diffeomorphism on $\mathbb{R}^n$. As a consequence, $\Phi$ is surjective.
\end{proposition}
Proof. We first proceed to show that there exists some matrix valued function $A(y):\mathbb{R}^n\to\mathbb{R}^{n\times n}$, such that $\Phi(y)-\Phi(0)=A(y)y$,~ for all $y\in\mathbb{R}^n$. To this end, for any given $y$, we define a function $h:[0,1]\to\mathbb{R}^n$ of $t$ by
 $$h(t)\overset{\triangle}{=}\Phi(x), ~~~x=ty,$$ then by the chain rule, it follows that the derivative of $h$ with respect to $t$ is
$$\frac{d h}{dt}=\frac{\partial \Phi}{\partial x}\frac{\partial x}{\partial t}=\left[\frac{\partial \Phi}{\partial x}(x)\right]y, ~~~~ x=ty.$$
It follows that
\begin{align*}
&\Phi(y)-\Phi(0)=h(1)-h(0)=\int_0^1\frac{d h}{dt} dt\\
=&\int_0^1\Big[\frac{\partial \Phi}{\partial x}(ty)\Big]y dt=A(y)y,
\end{align*}
where
$A(y)=\int_0^1 \big[\frac{\partial \Phi}{\partial x}(ty)\big] dt$.

Next, we proceed to show that $\text{Sym}[A(y)]\geq  R,~\forall y\in\mathbb{R}^n.$
Let $y\in\mathbb{R}^n$ be any given vector, then for any $\xi \in\mathbb{R}^n$, it is easy to see that
 \begin{align*}
 \xi^{\mathsf{T}}\text{Sym}[A(y)]\xi=& \xi^{\mathsf{T}}\text{Sym}\left[ \int_0^1\frac{\partial \Phi}{\partial x}(ty) dt\right] \xi\\
 =&\xi^{\mathsf{T}}\left\{\int_0^1 \text{Sym}\left[\frac{\partial \Phi}{\partial x}(ty)\right]dt\right\} \xi\\
 =&\int_0^1 \xi^{\mathsf{T}}\text{Sym} \left[\frac{\partial \Phi}{\partial x}(ty)\right]\xi dt
 \\ \geq &\int_0^1 \xi^{\mathsf{T}}R\xi dt
 =\xi^{\mathsf{T}}R\xi,
 \end{align*}
where the last but one relationship holds because by our assumption we have
$${\rm Sym}\left[\frac{\partial \Phi}{\partial x}(ty)\right]\geq R, ~\text{for all}~ t\in [0,1].$$
Finally, since $\xi\in\mathbb{R}^n$ is arbitrary, we conclude that $\text{Sym}[A(y)]\geq R$.
From this, we can obtain
 \begin{align}\label{norm}
 &\left|y^{\mathsf{T}}(\Phi(y)-\Phi(0))\right|=\left|y^{\mathsf{T}} A(y)y\right|\nonumber\\
 =&\left|y^{\mathsf{T}} \text{Sym}[A(y)]y\right| \geq  \lambda_{\min}(R) \|y\|^2,
 ~~\forall y\in\mathbb{R}^n,\end{align}
 where $\lambda_{\min}(R)>0$ is the smallest eigenvalue of $R$.
Then, it follows from (\ref{norm}) that $\|\Phi(y)-\Phi(0)\|\geq  \lambda_{\min}(R) \|y\|$, which in turn implies  $\lim_{\|y\|\to\infty} \|\Phi(y)\|=\infty.$ Combining this with the fact that $\text{Sym}\big[\frac{\partial \Phi}{\partial x}]$ is positive definite(and thus $\frac{\partial \Phi}{\partial x}$ is nonsingular),  we conclude that $\Phi$ is a global diffeomorphism on $\mathbb{R}^n$ by  Theorem A1. \hfill$\square$

\begin{proposition} Assume that $(k_i,k_p,k_d)\in \Omega_{pid}$, then the following $3n\times 3n$ matrix is positive definite, i.e.,
 \begin{align}\label{P}P\overset{\triangle}{=}
\begin{bmatrix}
2k_ik_p\underline bI_n&2k_ik_d\underline bI_n&k_iI_n\\
2k_ik_d\underline bI_n&~(2k_pk_d\underline b-k_i)I_n~&k_pI_n\\
k_iI_n&k_pI_n&k_dI_n
\end{bmatrix}>0.\end{align}
 Moreover, if we denote
 $$A\overset{\triangle}{=}\begin{bmatrix}
 0_n&I_n&0_n\\0_n&0_n&I_n\\-k_i\theta~ &~a-k_p\theta~&~b-k_d\theta\end{bmatrix},$$
where $0_n$ is the $n\times n$ zero matrix and $a$,  $b$ and $\theta$ are $n\times n$ constant matrices satisfying
 \begin{align}\label{common}\|a\|\le L_1, ~\|b\|\le L_2,~~{\rm Sym}[\theta] \geq \underline bI_n >0.\end{align}
Then there exists $\alpha>0$(depends on $(k_p,k_i,k_d,L_1,L_2,\underline b)$ only), such that   $PA+A^{\mathsf{T}}P\le -\alpha I_{3n}$ holds, where $I_{3n}$ is the $3n\times 3n$ identity matrix.
\end{proposition}

Proposition 4.2 implies that  the function $V(z)=z^{\mathsf{T}}Pz,z\in\mathbb{R}^{3n}$ with $P$ defined by (\ref{P}) is actually a common Lyapunov function of $\dot z=Az$ for all $a,b$ and $\theta$ satisfying (\ref{common}) in the matrix $A$. This fact will be used in the proof of Theorem 3.1 to be given shortly.

\begin{proposition} Assume that $(k_p,k_d)\in \Omega_{pd}$, then the following $2n\times 2n$ matrix  is positive definite, i.e.,
\begin{align}\label{P12}P\overset{\triangle}{=}\begin{bmatrix}
 2k_pk_d\underline bI_n~&~k_pI_n\\k_pI_n&k_dI_n
\end{bmatrix}>0.\end{align}
Moreover, if we denote
$$A\overset{\triangle}{=}\begin{bmatrix}0_n&I_n\\ a-k_p\theta~&~b-k_d\theta\end{bmatrix},$$ where $a$, $b$ and $\theta$ are $n\times n$ constant matrices satisfying
 $$\|a\|\le L_1, ~\|b\|\le L_2,~~{\rm Sym}[\theta] \geq \underline bI_n >0.$$
 Then there exists $\beta>0$(depends on $(k_p,k_d,L_1,L_2,\underline b)$ only), such that $PA+A^{\mathsf{T}}P\le -\beta I_{2n}$ holds, where $I_{2n}$ is the $2n\times 2n$ identity matrix.
\end{proposition}
\begin{proposition} Assume that $(k_p,k_i)\in \Omega_{pi}$, then the following $2n\times 2n$ matrix  is positive definite, i.e.,
\begin{align}\label{P1}P\overset{\triangle}{=}\begin{bmatrix}
 2k_pk_i\underline bI_n&k_iI_n\\k_iI_n&k_pI_n
\end{bmatrix}>0.\end{align}
Moreover, if we denote $$A\overset{\triangle}{=}\begin{bmatrix}0_n&I_n\\-k_i\theta~&~a-k_p\theta\end{bmatrix},$$ where $a$ and $\theta$ are $n\times n$ matrices satisfying $\|a\|\le L$ and ${\rm Sym}[\theta]\geq \underline bI_n >0$. Then there exists a constant $\gamma>0$(depends on $(k_i,k_p,L,\underline b)$ only), such that  $PA+A^{\mathsf{T}}P\le -\gamma I_{2n}$ holds.
\end{proposition}
The proofs of Propositions 4.2$-$4.4 are given in Section 6.1.

\textbf{Proof of Theorem 3.1}

Let $y^*\in \mathbb{R}^n$ be given, and suppose that the nonlinear function $f(x_1,x_2,u)\in \mathcal{F}_{L_1,L_2,\underline{b}}$. We first define a vector-valued function $\Phi:\mathbb{R}^n\to\mathbb{R}^n$ as follows: $$\Phi(u)\overset{\triangle}{=}f(y^*,0,u), ~~u\in\mathbb{R}^n.$$
Note that
 $$\text{Sym}\Big[\frac{\partial f}{\partial u}\Big]\geq \underline{b}I_n>0,~\forall x_1,x_2,u\in \mathbb{R}^n,$$ therefore $\text{Sym}[\frac{\partial \Phi}{\partial u}]\geq \underline{b}I_n>0,~\forall u\in \mathbb{R}^n$. By Proposition 4.1, we know that $\Phi$ is a global  diffeomorphism on $\mathbb{R}^n$. Hence, there exists a unique $u^*\in\mathbb{R}^n$, such that $\Phi(u^*)=f(y^*,0,u^*)=0$.

Now, suppose  that $(k_p,k_i,k_d)\in \Omega_{pid}$.  We first introduce the following notations:

$~~z_0(t)=\int_{0}^t e(s)ds-u^*/k_i,~~z_1(t)=e(t),~~z_2(t)=\dot e(t),$

then the PID controller (2) can be rewritten as
$$u(t)=k_iz_0(t)+k_pz_1(t)+k_dz_2(t)+u^*.$$
Define a function $g:\mathbb{R}^{3n}\to\mathbb{R}^n$ as follows:
\begin{align}\label{notation}
g(x_1,x_2,u)=&-f(y^*-x_1,-x_2,u+u^*).\end{align}
Notice that $f\in \mathcal{F}_{L_1,L_2,\underline{b}}$, it is easy to see that
\begin{equation}\label{ggg}
\big\|\frac{\partial{g}}{\partial{x_i}}\big\|\le L_i;~
\text{Sym}\Big[\frac{\partial g}{\partial u}\Big]\le -\underline{b}I_n,~\forall x_1,x_2,u\in\mathbb{R}^n.\end{equation}
Note that $\dot x_2=f(x_1,x_2,u)$, it follows that(we omit the variable $t$ for simplicity)
\begin{align*}
\dot {z}_2=&-\dot x_2=-f(x_1,x_2,u)\\=&-f(y^*-z_1,-z_2,k_iz_0+k_pz_1+k_dz_2+u^*)\\
=&g(z_1,z_2,k_iz_0+k_pz_1+k_dz_2).\end{align*}
Hence, we obtain
\begin{align}\label{g0}
\begin{cases}
\dot {z}_0=z_{1}, \\
\dot {z}_1=z_{2}, \\
\dot {z}_2=g(z_1,z_2,k_iz_0+k_pz_1+k_dz_2).
\end{cases}
\end{align}
Recall that $f(y^*,0,u^*)=0,$ we have $$g(0,0,0)=-f(y^*,0,u^*)=0.$$
Next, we proceed to show that there exists  a continuous function $\theta(z_0,z_1,z_2):\mathbb{R}^{3n}\to \mathbb{R}^{n\times n}$  such that
\begin{align}\label{abc}
&g(z_1,z_2,k_iz_0+k_pz_1+k_dz_2)\nonumber\\=&g(z_1,z_2,0)-\theta(z_0,z_1,z_2) (k_iz_0+k_pz_1+k_dz_2),
\end{align}
with  $~\text{Sym}[\theta(z_0,z_1,z_2)]\geq \underline{b}I_n>0,~~\forall (z_0,z_1,z_2)\in\mathbb{R}^{3n}.$

To this end, for fixed $z_0,z_1$ and $z_2$, we define a function $h(\tau)$ of $\tau \in [0,1]$ as follows: $$h(\tau)\overset{\triangle}{=}g(z_1,z_2,u),~~u=\tau(k_iz_0+k_pz_1+k_dz_2).$$ By the chain rule, we have
$$ \frac{d h}{d\tau}=\frac{\partial g}{\partial u}\frac{\partial u}{\partial \tau}
 \!=\!\left[\frac{\partial g}{\partial u}\!\big(z_1,\!z_2,\!u\big)\right](k_iz_0+k_pz_1+k_dz_2).$$
Then it follows that
\begin{align*}
 &g(z_1,z_2,k_iz_0+k_pz_1+k_dz_2)-g(z_1,z_2,0)
 \\=&h(1)-h(0)=\int_0^1 \frac{d h}{d\tau}d\tau\\
 =&\left[\int_0^1 \frac{\partial g}{\partial u}\big(z_1,z_2,u\big)d\tau\right](k_iz_0+k_pz_1+k_dz_2).
 \end{align*}
 Therefore, if we denote $\theta$ as follows: \begin{align*}\theta=\int_0^1 -\frac{\partial g}{\partial u}(z_1,z_2,u)d\tau,~~u=\tau(k_iz_0+k_pz_1+k_dz_2),\end{align*}
 then it is easy to see that $\theta=\theta(z_0,z_1,z_2)$ is a continuous function of $z_0,z_1$ and $z_2$. Moreover, (\ref{abc}) is satisfied obviously.
 Note that $\text{Sym}[-\frac{\partial g}{\partial u}]\geq  \underline{b}I_n$ and the integral expression of $\theta$, similar to the proof of Proposition 4.1, we conclude that
  $$\text{Sym}[\theta(z_0,z_1,z_2)]\geq \underline{b}I_n>0,~~\forall (z_0,z_1,z_2)\in\mathbb{R}^{3n}.$$
Similarly, recall that $g(0,0,0)=0$, we can express $g(z_1,z_2,0)$ as follows:
\begin{align}\label{ac}
g(z_1,z_2,0)=&\left[g(z_1,0,0)-g(0,0,0)\right]\nonumber\\
&+\left[g(z_1,z_2,0)-g(z_1,0,0)\right]\nonumber\\=&~a(z_1) z_1+b(z_1,z_2)z_2,
\end{align}
 where $a(z_1)$ and $b(z_1,z_2)$ are $n\times n$ matrices defined by
 $$a(z_1)=\int_0^1 \frac{\partial g}{\partial x_1}(\tau z_1,0,0)d\tau,$$ and $$b(z_1,z_2)=\int_0^1 \frac{\partial g}{\partial x_2}(z_1,\tau z_2,0)d\tau.$$
 Hence,  by (\ref{ggg}), it is easy to obtain
    $$\|a(z_1)\|\le \int_0^1\| \frac{\partial g}{\partial x_1}(\tau z_1,0,0)\|d\tau\le L_1,~~\forall z_1\in\mathbb{R}^n,$$
and $$\|b(z_1,z_2)\|\le\! \int_0^1 \!\|\frac{\partial g}{\partial x_2}(z_1,\tau z_2,0)\|d\tau \le L_2,~\forall z_1,z_2\in\mathbb{R}^n.$$
By (\ref{abc})-(\ref{ac}),  equation (\ref{g0}) turns into
\begin{equation}\label{g1}
\begin{cases}
\dot {z}_0=&z_1\\
\dot {z}_1=&z_2\\
\dot {z}_2=&a(z_1)z_1+b(z_1,z_2)z_2\\&-\theta(z_0,z_1,z_2) (k_iz_0+k_pz_1+k_dz_2).\\
\end{cases}
\end{equation}

For simplicity, let us  denote
$a=a(z_1)$, $b=b(z_1,z_2)$, $\theta=\theta(z_0,z_1,z_2)$ and
$z=(z_0^{\mathsf{T}},z_1^{\mathsf{T}},z_2^{\mathsf{T}})^{\mathsf{T}}$ and define  $$A(z)\overset{\triangle}{=}\begin{bmatrix}0_n&I_n&0_n\\0_n&0_n&I_n\\-k_i\theta ~&a-k_p\theta~&b-k_d\theta\end{bmatrix}.$$
Then (\ref{g1}) can be rewritten as $\dot z=A(z) z, ~z\in\mathbb{R}^{3n}$.

Now, we construct the following Lyapunov function:
$$V(z)=z^{\mathsf{T}}Pz, ~~z=(z_0^{\mathsf{T}},z_1^{\mathsf{T}},z_2^{\mathsf{T}})^{\mathsf{T}},$$ where $P$ is a $3n\times 3n$ constant matrix defined by (\ref{P}).
Then, the time derivative of $V(\cdot) $ along the trajectories of (\ref{g1}) is given by
$$\dot V(z(t))=z^{\mathsf{T}}(t)\left(PA(z(t))+A^{\mathsf{T}}(z(t))P\right)z(t).$$
By Proposition 4.2, we know that $P>0$. Moreover, recall that the $n\times n$ matrices $a(z_1)$, $b(z_1,z_2)$ and $\theta(z_0,z_1,z_2)$  satisfies
 $$\|a(z_1)\|\le L_1,~\|b(z_1,z_2)\|\le L_2,~\text{Sym}[\theta(z_0,z_1,z_2)]\geq \underline{b}I_n,$$for all $z_0$, $z_1$ and $z_2$. Thus, by Proposition 4.2 again, we know that there exists a time-invariant constant $\alpha>0$ which depends on $(k_p,k_i,k_d,L_1,L_2,\underline b)$ only,  such that  $$PA(z)+A^{\mathsf{T}}(z)P\le -\alpha I_{3n},~~ \text{for all}~~ z\in\mathbb{R}^{3n}.$$
As a consequence, we have  $\dot V(z(t))\le -\alpha \|z(t)\|^2$, which implies the origin $0\in\mathbb{R}^{3n}$ is globally exponentially stable of (\ref{g1}).
In fact, it is easy to obtain
$$\dot V(z(t))\le -\alpha \|z(t)\|^2\le -\alpha V(z(t))/\lambda_{\max}(P),$$ which in turn gives
\begin{align*}\lambda_{\min }(P)\|z(t)\|^2\le& V(z(t))\le e^{-\frac{\alpha t}{\lambda_{\max}(P)}}V(z(0))\\
\le&\lambda_{\max }(P) e^{-\frac{\alpha t}{\lambda_{\max}(P)}}\|z(0)\|^2.\end{align*}
As a consequence, we have
$$\|z(t)\|\le\sqrt{\lambda_{\max }(P)/\lambda_{\min }(P) }e^{-\frac{\alpha t}{2\lambda_{\max}(P)}}\|z(0)\|.$$
Notice that $z(t)=(z^{\mathsf{T}}_0(t),z^{\mathsf{T}}_1(t),z^{\mathsf{T}}_2(t))^{\mathsf{T}}$, and that $z_1(t)=e(t)$, $z_2(t)=\dot e(t)$ and $z_0(0)=-u^*/k_i$, it is easy to obtain
\begin{align}\label{et}
&\|e(t)\|+\|\dot e(t)\|\le
\sqrt 2 \|z(t)\|\nonumber\\\le& M_1e^{-\lambda t} \big(\|e(0)\|+\|\dot e(0)\|+\|u^*\|/k_i\big),
\end{align}
where $M_1=\sqrt{2\lambda_{\max }(P)/\lambda_{\min }(P) }$ and $\lambda=\alpha /(2\lambda_{\max}(P))$. Denote
$$M=\max\{M_1,M_1/k_i\},$$ from (\ref{et}), we obtain
\begin{align*}
&\|e(t)\|+\|\dot e(t)\|\nonumber\\\le&  Me^{-\lambda t} \big(\|e(0)\|+\|\dot e(0)\|+\|u^*\|\big), ~\text{for all}~t\geq 0,
\end{align*}
 where $M$ is a constant depends on  $(k_p,k_i,k_d,\underline b)$ only and $\lambda$ is a constant depends on  $k_p,k_i,k_d,L_1,L_2$ and $\underline b$ only. Hence,  the proof of Theorem 3.1 is completed. \hfill$\square$

\textbf{Proof of Theorem 3.2}

Suppose that $f\in \mathcal{F}_{L_1,L_2,\underline{b},y^*}$ and $(k_p,k_d)\in \Omega_{pd}$.  Denote $z_1(t)=e(t)=y^*-x_1(t),~z_2(t)=\dot e(t)=-x_2(t)$ and $~g(x_1,x_2,u)=-f(y^*-x_1,-x_2,u),$ then the PD controlled system (1) and (5) turns into
\begin{equation}\label{g000}
\begin{cases}
\dot {z}_1=z_2\\
\dot {z}_2=g(z_1,z_2,k_pz_1+k_dz_2).\\
\end{cases}
\end{equation}
From $f(y^*,0,0)=0,$ we have $g(0,0,0)=-f(y^*,0,0)=0,$
which implies $0\in\mathbb{R}^{2n}$ is an equilibrium of (\ref{g000}).

Next, similar to the proof of Theorem 3.1, we know that $g(z_1,z_2,k_pz_1+k_dz_2)$ can be expressed by
 \begin{align*}&g(z_1,z_2,k_pz_1+k_dz_2)\\
 =&a(z_1) z_1+b(z_1,z_2)z_2-\theta(z_1,z_2) (k_pz_1+k_dz_2),\end{align*} where
  $a(z_1)$, $b(z_1,z_2)$ and $\theta(z_1,z_2)$ are $n\times n$ matrices  satisfying
 \begin{equation}\label{gg}
 \|a(z_1)\|\le L_1,~~\|b(z_1,z_2)\|\le L_2,~~\text{Sym}[\theta(z_1,z_2)]\geq \underline{b}I_n,
 \end{equation}
for all $z_1,z_2\in\mathbb{R}^{n}$.
Therefore,  (\ref{g000}) turns into
\begin{equation}\label{g111}
\begin{cases}
\dot {z}_1=z_2\\
\dot {z}_2=a(z_1)z_1+b(z_1,z_2)z_2-\theta(z_1,z_2) (k_pz_1+k_dz_2)\\
\end{cases}
\end{equation}
Denote $z=(z_1^{\mathsf{T}},z_2^{\mathsf{T}})^{\mathsf{T}}$ and $$A(z)\overset{\triangle}{=}\begin{bmatrix}0_n&I_n\\a(z_1)\!-\!k_p\theta(z_1,z_2)&~b(z_1,z_2)\!-\!k_d\theta(z_1,z_2)\end{bmatrix},$$
then (\ref{g111}) can be rewritten as $\dot z=A(z) z$.

Now, we construct the following Lyapunov function:
$$V(z)=z^{\mathsf{T}}Pz,$$ where  $P$ is a $2n\times 2n$ constant matrix defined by (\ref{P12}).
Therefore, the time
derivative of $V $ along the trajectories of (\ref{g111}), is given by $\dot V(z)=z^{\mathsf{T}}(PA(z)+A(z)^{\mathsf{T}}P)z$.
By Proposition 4.3, we know that $P>0$ and there exists $\beta>0$ such that  $PA(z)+A^{\mathsf{T}}(z)P\le -\beta I_{2n}$, for all $z\in\mathbb{R}^{2n}$.
Similarly, we have
$$\|z(t)\|\le\sqrt{\lambda_{\max }(P)/\lambda_{\min }(P) }e^{-\frac{\beta t}{2\lambda_{\max}(P)}}\|z(0)\|.$$
Notice that $z(t)=(z^{\mathsf{T}}_1(t),z^{\mathsf{T}}_2(t))^{\mathsf{T}}$ and that $z_1(t)=e(t)$, and $z_2(t)=\dot e(t)$, it is easy to obtain
\begin{align*}
\|e(t)\|+\|\dot e(t)\|\le&\sqrt 2\|z(t)\|\le  Me^{-\lambda t} \big(\|e(0)\|+\|\dot e(0)\|\big),
\end{align*}
where
$$M=\sqrt{\frac{2\lambda_{\max }(P)}{\lambda_{\min }(P)}},~~~\lambda=\frac{\beta}{2\lambda_{\max}(P)}.$$
Hence, the proof of Theorem 3.2 is completed.  \hfill$\square$

\textbf{Proof of Theorem 3.3}

Let $y^*\in \mathbb{R}^n$ be  given and suppose that $f(x,u)\in \mathcal{F}_{L,\underline{b}}$ and $(k_p,k_i)\in \Omega_{pi}$.  Define a map $\Phi:\mathbb{R}^n\to\mathbb{R}^n$ as follows: $\Phi(u)\overset{\triangle}{=}f(y^*,u), ~u\in\mathbb{R}^n.$
Note that $$\text{Sym}\Big[\frac{\partial \Phi}{\partial u}\Big]=\text{Sym}\Big[\frac{\partial f}{\partial u}(y^*,u)\Big]\geq \underline{b}I_n,~\forall u\in \mathbb{R}^n,$$ it follows from Proposition 4.1 that there exists a unique $u^*\in\mathbb{R}^n$, such that  $f(y^*,u^*)=0$.
 Denote \begin{align}\label{notation1}
z_0(t)=&\int_{0}^t e(s)ds-\frac{u^*}{k_i},~z_1(t)=e(t)=y^*-x(t),\nonumber
\\~g(x,u)&=-f(y^*-x,u+u^*),\end{align} then  $g(0,0)=-f(y^*,u^*)=0$ and it is not difficult to obtain
\begin{equation}\label{g00}
\begin{cases}
\dot {z}_0=z_1\\
\dot {z}_1=g(z_1,k_iz_0+k_pz_1).\\
\end{cases}
\end{equation}
Notice that $f\in \mathcal{F}_{L,\underline{b}}$ and $g(x,u)=-f(y^*-x,u+u^*)$, then it is easy to see that
\begin{equation}\label{g}\big\|\frac{\partial{g}}{\partial{x}}\big\|\le L,~~~
\text{Sym}\Big[\frac{\partial g}{\partial u}\Big]\le -\underline{b}I_n<0,~~~~\forall x,u\in\mathbb{R}^n.
\end{equation}

Next, similar to the proof of Theorem 3.1, there exist  continuous functions $a(z_1):\mathbb{R}^{n}\to \mathbb{R}^{n\times n}$ and $\theta(z_0,z_1):\mathbb{R}^{2n}\to \mathbb{R}^{n\times n}$  such that
\begin{equation}\label{ab}
g(z_1,k_iz_0+k_pz_1)=a(z_1)z_1-\theta(z_0,z_1) (k_iz_0+k_pz_1),
\end{equation}
with  $$\|a(z_1)\|\le L,~~~\text{Sym}[\theta(z_0,z_1)]\geq \underline{b}I_n>0,~~~\forall z_0,z_1\in\mathbb{R}^{n}.$$
Thus,  if we denote $a=a(z_1)$, $\theta=\theta(z_0,z_1)$ and $z=(z_0^{\mathsf{T}},z_1^{\mathsf{T}})^{\mathsf{T}}$ and  $$A(z)\overset{\triangle}{=}\begin{bmatrix}0_n&I_n\\-k_i\theta(z_0,z_1) ~&~a(z_1)-k_p\theta(z_0,z_1)\end{bmatrix}.$$
Then (\ref{g00}) can be rewritten as $\dot z=A(z) z, ~z\in\mathbb{R}^{2n}$.

Now, we construct the following Lyapunov function:
$$V(z)=z^{\mathsf{T}}Pz,$$ where  $P$ is a $2n\times 2n$ constant matrix defined by (\ref{P1}).
Then, the time derivative of $V $ along the trajectories of (\ref{g00}) is given by
$$\dot V(z(t))=z^{\mathsf{T}}(t)\left(PA(z(t))+A^{\mathsf{T}}(z(t))P\right)z(t).$$
By Proposition 4.4, we know that $P$ is positive definite and there exists $\gamma>0$ such that  $PA(z)+A^{\mathsf{T}}(z)P\le -\gamma I_{2n}$, for all $z\in\mathbb{R}^{2n}$.
By the comparison lemma, it is easy to obtain
$$\|z(t)\|\le\sqrt{\lambda_{\max }(P)/\lambda_{\min }(P) }e^{-\frac{\gamma t}{2\lambda_{\max}(P)}}\|z(0)\|.$$
Notice that $z(t)=(z_0(t),z_1(t))$ and that $z_1(t)=y^*-x(t)$, and $z_0(0)=u^*$, we have
\begin{align*}
\|x(t)-y^*\|\le\|z(t)\|\le Me^{-\lambda t} \left(\|x(0)-y^*\|+\|u^*\|\right),
\end{align*}
where $M=\sqrt{\lambda_{\max }(P)/\lambda_{\min }(P) }$ and $\lambda=\gamma /(2\lambda_{\max}(P))$.

Next, we proceed to show that if the dimension of $x$ and $u$ are both one, then the PI parameters $(k_p,k_i)$ can be chosen from the larger and necessary set $$\Omega'_{pi}=\left\{( k_p,k_i)\in \mathbb{R}^2\big| k_p\underline b>L,~~k_i>0\right\}.$$

\textbf{Sufficiency:} Let $y^*\in\mathbb{R}$, $f\in \mathcal{F}_{L,\underline b}$ and $(k_p,k_i)\in \Omega'_{pi}$.

Since $\frac{\partial f}{\partial u}(x,u)\geq \underline{b}>0,~\forall x,u\in \mathbb{R}$, it follows that  $$\lim_{u\to \infty}f(y^*,u)=\infty, ~~\lim_{u\to -\infty}f(y^*,u)=-\infty.$$ By the intermediate value theorem,   there exists  a unique $u^*\in\mathbb{R}$, such that $f(y^*,u^*)=0$.

Denote $z_0(t)=\int_0^t e(s)ds -u^*/k_i$, $z_1(t)=e(t)=y^*-x(t)$ and $g(x,u)=-f(y^*-x,u+u^*)$,   then (\ref{g00}) is satisfied with
$$|\frac{\partial g}{\partial x}|\le L,
~~\frac{\partial g}{\partial u}\le -\underline b<0, ~~\forall x,u\in\mathbb{R}.$$
Since $g(0,0)=-f(y^*,u^*)=0$, we know that $(0,0)$ is an equilibrium of (\ref{g00}).

We proceed to apply the Markus-Yamabe's theorem(Theorem A2) to show  $(0,0)$ is a globally asymptotically stable equilibrium of (\ref{g00}).
Denote   the vector field of $(\ref{g00})$ by $G(z_0,z_1)$, i.e.,
\begin{align*}
G(z_0,z_1)=\begin{bmatrix}
z_1\\g(z_1,k_iz_0+k_pz_1)
\end{bmatrix},~~~(z_0,z_1)\in\mathbb{R}^2.
\end{align*} Then the Jacobian matrix of $G$ is
$$DG(z_0,z_1)=\begin{bmatrix}
 ~~0& ~~1\\ k_i\dfrac{\partial{g}}{\partial{u}}~~&
 ~~\dfrac{\partial{g}}{\partial{x}}+k_p\dfrac{\partial g}{\partial u}\end{bmatrix}.$$
Since $(k_p,k_i)\in \Omega'_{pi}$, and by $|\frac{\partial{g}}{\partial{x}}|\le L,
\frac{\partial g}{\partial u}\le -\underline{b}<0$, we can easily see that
the trace of $DG(z_0,z_1)$ satisfies
 \begin{align*}
 \frac{\partial{g}}{\partial{x}}+k_p\dfrac{\partial g}{\partial u}\le L-k_p\underline b<0,~~(z_0,z_1)\in \mathbb{R}^2,
 \end{align*}
 and the determinant of $DG(z_0,z_1)$ satisfies
  \begin{align*}
 -k_i\dfrac{\partial{g}}{\partial{u}}>0,~~~~(z_0,z_1)\in \mathbb{R}^2,
  \end{align*}
   which implies that the two eigenvalues of $DG(z_0,z_1)$ have negative real parts for any $(z_0,z_1)$. By Theorem A2, we conclude that $(0,0)$ is a globally  asymptotically stable equilibrium of (\ref{g00}), thus $\lim_{t \to \infty} x(t)=y^{*}$ for any initial value $x(0)\in \mathbb{R}$.

\textbf{Necessity:}   We use  contradiction argument and  assume that $e(t)\rightarrow 0$ for all $f\in \mathcal{F}_{L,\underline b}$  and for all setpoint $y^*\in\mathbb{R}$  but $(k_{p},k_{i})\notin \Omega'_{pi}$. We  consider two cases separately.

(i) If $k_i=0$, let $f(x,u)=Lx+\underline b u\in \mathcal{F}_{L,\underline b}$, then $$\dot{e}(t)=-\dot x(t)=(L-\underline bk_p)e(t)-Ly^*.$$ Obviously, $e(t)$ does not converge to $0$ as long as $y^*\neq 0$ for any initial value $e(0)\in \mathbb{R}$.

(ii) If $k_i\neq 0$, we take $f(x)=Lx+\underline b u$ and $y^*=0$. If we denote $e_0(t)=\int_0^t e(s)ds$, then we have
\begin{equation}
\begin{cases}
\dot{e}_0=e\\
\dot{e}=-k_i\underline b e_0+(L-k_p\underline b) e.
\end{cases}
\end{equation}
  Since $(k_{p},k_{i})\notin \Omega'_{pi}$, then there exists at least one eigenvalue  $\lambda$ of $A=\begin{bmatrix}0&1\\-k_i\underline b&L-k_p\underline b\end{bmatrix}$, whose real part satisfies $\Re(\lambda)\geq 0$,
 which contradicts with our assumption that $e(t)\rightarrow 0$ for all initial states $x(0)\in\mathbb{R}$.  \hfill$\square$

\textbf{Remark 4.1.}  From the proof of Theorems 3.1-3.3, one can see that some seemly linearity-like methods are used in this paper, which is not surprising since the classical PID control is a linear feedback. What somewhat surprising to us are the following two facts: a) the linear PID has proven to have the ability to globally stabilize and regulate a wide class of nonlinear systems; b) some  seemly linearity-like methods can still be used to analyze the  PID controlled strongly coupled non-affine nonlinear uncertain systems. Moreover, the linearity-like methods used in Section 4 cannot be directly borrowed from the conventional linear system theory for at least the following two reasons: 1) the construction of the Lyapunov functions in Theorems 3.1$-$3.3 is far from obvious because the sophisticated  structure of the positive definite matrices $P$ in (9), (11)-(12), which are not the conventionally assumed solutions of certain matrix equations or inequalities; 2) The proofs of Theorems 3.1 and 3.3 need a basic topological result known as the Hadamard global inverse function theorem, and the proof of Theorem 3.3 also needs the Markus-Yamabe theorem which had been a longstanding problem in differential equations and has been rarely used in control theory before.

\section{Conclusion}
In this paper, we have presented a mathematical investigation on the foundation of the classical linear PID (PD or PI) controller for several classes of MIMO non-affine uncertain nonlinear dynamical systems,  found a three (two) dimensional
set within which the PID (PD or PI) parameters can be chosen arbitrarily to stabilize and regulate the system globally. Moreover, a simple necessary and sufficient condition for the choice of the PI controller parameters is given for a class of one dimensional non-affine uncertain systems.  The PID design rules given in this paper are quite simple, which shows that the PID parameters are not necessarily to be of high gain. Our results have improved significantly the existing design methods and theoretical results in the related literature. Of course,  many interesting problems
still remain open. For example, it would be natural to ask whether the extended PID controller discussed in \cite{zhao2020} can globally stabilize more general non-affine uncertain systems under the similar assumptions on the system unknown functions, and how  to further ``optimize'' the PID parameters from  $\Omega_{pid}$ to get better transient performance.  It would also be interesting to consider more complicated  situations such as saturation, dead-zone, time-delayed inputs, sampled-data PID controllers under a prescribed sampling rate, etc.
\section{Appendix}

\subsection{Proof of Propositions}
To prove Proposition 4.2, we need the following three lemmas.
\begin{lemma}
Suppose $D=D^{\mathsf{T}}\in\mathbb{R}^{m\times m}$ and $E=E^{\mathsf{T}}\in\mathbb{R}^{n\times n}$ are two invertible matrices and $B\in\mathbb{R}^{m\times n}$, then $\begin{bmatrix}D&B\\B^{\mathsf{T}}&E\end{bmatrix}>0$ if and only if
$D>0$ and $E-B^{\mathsf{T}}D^{-1}B>0$.
\end{lemma}

Lemma 6.1 can be found  in \cite{horn}( see Theorem 7.7.6). By Lemma 6.1, we can easily  obtain the following results:
\begin{lemma}
Suppose  $D\in\mathbb{R}^{m\times m}$ and $E\in\mathbb{R}^{n\times n}$ are two positive definite matrices and  $B\in\mathbb{R}^{m\times n}$. If
 $\lambda_{\min}(D)\lambda_{\min}(E)>\|B\|^2$, where $\lambda_{\min}(D)$ denotes the smallest eigenvalue of $D$, then $\begin{bmatrix}D&B\\B^{\mathsf{T}}&E\end{bmatrix}>0$.
\end{lemma}
Proof.  From Lemma 6.1, it suffices to show that $E-B^{\mathsf{T}}D^{-1}B>0$. First, it is easy to see that
$\lambda_{\min} (D)D^{-1}\le I_m$. Let $x\in\mathbb{R}^n$,
then we can obtain
\begin{align}\label{xe}
 \lambda_{\min} (E)\|x\|^2\le& x^{\mathsf{T}}Ex,\nonumber\\
x^{\mathsf{T}}B^{\mathsf{T}}D^{-1}Bx\le&\frac{\|Bx\|^2}{\lambda_{\min} (D)}\le \frac{\|B\|^2}{\lambda_{\min} (D)}\|x\|^2.
\end{align}
From (\ref{xe}), we have $$x^{\mathsf{T}}(E-B^{\mathsf{T}}D^{-1}B)x\geq \big(\lambda_{\min}(E)-\frac{\|B\|^2}{\lambda_{\min} (D)}\big)\|x\|^2.$$ Thus, $E-B^{\mathsf{T}}D^{-1}B>0$.  \hfill$\square$

\begin{lemma}
Let $F:K\to \mathbb{R}^{n\times n}$ be a continuous (matrix valued) function, where $K\subset\mathbb{R}^m$ is a compact set. Assume that for any $ z \in K$, $F(z)$ is positive definite. Then  $$\inf_{z\in K}\lambda_{\min} \big(F(z)\big)>0,$$
where $\lambda_{\min}(F)$ denotes the smallest eigenvalue of $F$.
\end{lemma}
Proof. For a positive definite matrix $F$, it is well-known that  $\lambda_{\max}(F^{-1})=1/\lambda_{\min}(F)$ and $\|F^{-1}\|=\lambda_{\max}(F^{-1})$, where $\lambda_{\max}(F^{-1})$ denotes the largest eigenvalue of $F^{-1}$. Therefore \begin{equation}\label{F}\lambda_{\min}(F)=1/\|F^{-1}\|,~~F>0.\end{equation}
From (\ref{F}), we can see that to prove $ \inf_{z\in K}\lambda_{\min} F(z)>0$, it suffices to show $\sup_{z\in K}\|F^{-1}(z)\|<\infty.$ This comes immediately since the function
$z\mapsto \|F^{-1}(z)\|$ is continuous for $z\in K$ and $K$ is a compact set in $\mathbb{R}^m$. \hfill$\square$

 \textbf{Proof of Proposition 4.2}

 \textbf{Step 1.} We first show that the matrix $P$ defined by (\ref{P}) is positive definite.
 For this, it suffices to show that the following $3\times3$ matrix $P_0$ defined by
 \begin{align*}P_0\overset{\triangle}{=}\begin{bmatrix}
2k_ik_p\underline b&2k_ik_d\underline b&k_i\\
2k_ik_d\underline b&~~2k_pk_d\underline b-k_i~&k_p\\
k_i&k_p&k_d
\end{bmatrix}
\end{align*}
is positive definite,
since the matrix $P$ is also symmetric and shares the same spectrum with  $P_0$.
Moreover, to prove the positiveness of $P_0$, we need only to verify the following three inequalities:
$$2k_ik_p\underline b>0,~~\det \begin{bmatrix}
2k_ik_p\underline b&2k_ik_d\underline b\\
2k_ik_d\underline b&2k_pk_d\underline b\!-\!k_i
\end{bmatrix}>0,~~\det (P_0)>0.$$
Firstly, recall $k_i>0$, $k_p>0$ and $\underline b>0$, thus $2k_ik_p\underline b>0$.
Secondly, note that $k_p^2>2k_ik_d$ and $k_d^2>k_p/\underline b$, we know that
\begin{align*}
&~\det \begin{bmatrix}
2k_ik_p\underline b&~2k_ik_d\underline b\\
2k_ik_d\underline b&~2k_pk_d\underline b\!-\!k_i
\end{bmatrix}\\
=&2k_ik_p\underline b(2k_pk_d\underline b\!-\!k_i)-(2k_ik_d\underline b)^2\\
=&4k_ik_d(k_p^2-k_ik_d)\underline b^2-2k_i^2k_p\underline b\\
>& 4k_i^2k_d^2\underline b^2-2k_i^2k_p\underline b\\
=&2k_i^2\left(2k_d^2\underline b-k_p\right)\underline b
>0.
\end{align*}
Thirdly, by some simple calculations, we have
\begin{align}
\det (P_0)=k_i\left(4k_p^2k_d^2\underline b^2+k_i^2-2k_p^3\underline b-4k_ik_d^3\underline b^2\right).
\end{align}
Note that
$$\sqrt{k_p/\underline b}<k_d<\frac{k_p^2}{2k_i},$$ we define a function $\rho(x)$
as follows:
$$\rho(x)=4k_p^2\underline b^2x^2+k_i^2-2k_p^3\underline b-4k_i\underline b^2x^3, ~~\sqrt{k_p/\underline b}<x<\frac{k_p^2}{2k_i}.$$
Then, the derivative of $\rho(\cdot)$ satisfies
\begin{align*}
\frac{d\rho(x)}{dx}=&4x\left(2k_p^2-3k_ix\right)\underline b^2\\
>&4x\underline b^2\Big(2k_p^2-3k_i\frac{k_p^2}{2k_i}\Big)\\
>&2x\underline b^2k_p^2>0,~~\sqrt{k_p/\underline b}<x<\frac{k_p^2}{2k_i}.
\end{align*}
Hence, $\rho(x)$ is strictly increasing on the interval $(\sqrt{k_p/\underline b},\frac{k_p^2}{2k_i})$.
Consequently, by the monotonicity of $\rho(x)$, it can be deduced that
\begin{align*}
\frac{\det (P_0)}{k_i}=&\rho(k_d)>
  \rho \Big(\sqrt{k_p/\underline b}~\!\Big)\\
  =&2k_p\Big(k_p^2-2k_i\sqrt{k_p/\underline b}~\!\Big)\underline b+k_i^2\\
  >&2k_p\Big(2k_ik_d-2k_i\sqrt{k_p/\underline b}~\!\Big)\underline b+k_i^2\\
  >& 4k_ik_p\Big(k_d-\sqrt{k_p/\underline b}~\!\Big)\underline b+k_i^2 >0.\end{align*}
Thus, $P$ is a positive definite matrix.

\textbf{Step 2.}
Next, we calculate $PA+A^{\mathsf{T}}P$.

For simplicity, we denote
\begin{equation}
\hat \theta=\text{Sym}[\theta]=\frac{\theta+\theta^{\mathsf{T}}}{2},\end{equation}
then by some elementary calculations, it can be obtained that
\begin{align}\label{Q}
&Q\overset{\triangle}{=}-(PA+A^{\mathsf{T}}P)\nonumber\\=&\begin{bmatrix}
2k_i^2\hat \theta&~2k_ik_p\hat \theta-m_2&~2k_ik_d\hat \theta-m_3\\
2k_ik_p\hat \theta-m_2^{\mathsf{T}}&~2k_p^2\hat \theta-m_4&~2k_pk_d\hat \theta-m_5\\
2k_ik_d\hat \theta-m_3^{\mathsf{T}}&~2k_pk_d\hat \theta-m_5^{\mathsf{T}}&~ 2k_d^2\hat \theta-m_6
\end{bmatrix},
\end{align}
where $m_i$, $i=2,\cdots,6$ are $n\times n$ matrices defined by
 \begin{align} \label{c7}m_2=&2k_ik_p\underline bI_n+k_ia,\nonumber\\m_3=&2k_ik_d\underline bI_n+k_ib,\nonumber\\
  ~m_4=&4k_ik_d\underline bI_n+k_p \left(a+a^{\mathsf{T}}\right),\nonumber\\
 m_5=&2k_pk_d\underline b I_n+k_pb+k_da^{\mathsf{T}},~~~\nonumber\\  m_6=&2k_pI_n+k_d \left(b+b^{\mathsf{T}}\right).\end{align}
\textbf{Step 3.}
Denote $\theta_0=\underline b I_n$, we proceed to show that the matrix $Q$ defined by (\ref{Q}) satisfies the  inequality $Q\geq Q_0$, where
\begin{equation} \label{q0}Q_0\overset{\triangle}{=}
\begin{bmatrix}
2k_i^2 \theta_0&~2k_ik_p \theta_0-m_2&~2k_ik_d \theta_0-m_3\\
2k_ik_p \theta_0-m_2^{\mathsf{T}}&~2k_p^2 \theta_0-m_4&~2k_pk_d\theta_0-m_5\\
2k_ik_d \theta_0-m_3^{\mathsf{T}}&~2k_pk_d \theta_0-m_5^{\mathsf{T}}&~ 2k_d^2 \theta_0-m_6
\end{bmatrix}.\end{equation}

Notice that  \begin{align*}Q-Q_0
=\begin{bmatrix}
2k_i^2&&&2k_ik_p&&&2k_ik_d\\
2k_ik_p&&&2k_p^2&&&2k_pk_d\\
2k_ik_d&&&2k_pk_d&&&2k_d^2
\end{bmatrix}\otimes (\hat \theta-\theta_0),\end{align*}
where $\otimes$ denotes the  Kronecker product.
It can be verified  that
$$\begin{bmatrix}
2k_i^2&&&2k_ik_p&&&2k_ik_d\\
2k_ik_p&&&2k_p^2&&&2k_pk_d\\
2k_ik_d&&&2k_pk_d&&&2k_d^2
\end{bmatrix}\geq 0.$$  On the other hand, recall that $\hat\theta-\theta_0\geq 0$. Since the Kronecker product of two positive semi-definite matrices is also positive semi-definite, thus  $Q-Q_0\geq 0.$

\textbf{Step 4.} We proceed to show that $Q_0>0$ for all $\|a\|\le L_1$ and $\|b\|\le L_2$.

First, by the definition of $Q_0$ in (\ref{q0}), we know that
\begin{align}\label{q00}
Q_0=
\begin{bmatrix}
2k_i^2\underline bI_n&-k_ia&-k_ib\\
-k_ia^{\mathsf{T}}&2k_1I_n-2 k_p\hat a&~-(k_pb+k_da^{\mathsf{T}})\\
-k_ib^{\mathsf{T}}&-(k_pb^{\mathsf{T}}+k_da)&~2k_2 I_n-2k_d\hat b\end{bmatrix},
\end{align}
where
\begin{align*}
k_1\overset{\triangle}{=}&(k_p^2-2k_ik_d)\underline b, ~~k_2\overset{\triangle}{=}k_d^2\underline b-k_p,\\
\hat a\overset{\triangle}{=}&\left(a+a^{\mathsf{T}}\right)/2,~~ ~~ \hat b\overset{\triangle}{=}\left(b+b^{\mathsf{T}}\right)/2.
\end{align*}

Denote
\begin{align*}D\overset{\triangle}{=}&2k_i^2\underline bI_n,~~~ B\overset{\triangle}{=}[-k_i a,~-k_ib],\\
E\overset{\triangle}{=}&\begin{bmatrix}
2k_1I_n-2 k_p\hat a&~-(k_pb+k_da^{\mathsf{T}})\\
-(k_pb^{\mathsf{T}}+k_da)&~2k_2 I_n-2k_d\hat b\end{bmatrix},\end{align*}
then it follows from (\ref{q00}) that $Q_0=\begin{bmatrix}D&B\\B^{\mathsf{T}}&E\end{bmatrix}$.

Note that $D>0$, by Lemma 6.1, to prove  $Q_0>0$, it suffices to show that
 \begin{equation}\label{E}
 E-B^{\mathsf{T}}D^{-1}B>0.\end{equation}
By some calculations, it is not difficult to obtain
  \begin{align*}E-B^{\mathsf{T}}D^{-1}B\overset{\triangle}{=}&\begin{bmatrix}D_1&B_1\\B_1^{\mathsf{T}}&E_1\end{bmatrix}
,\end{align*}
where $D_1$, $B_1$ and $E_1$ have the following forms:
\begin{align*}
&D_1=2k_1I_n-2k_p\hat a-\frac{a^{\mathsf{T}}a}{2\underline b},\\
&B_1=-\Big(k_pb+k_da^{\mathsf{T}}+\frac{a^{\mathsf{T}}b}{2\underline b}\Big),\\
&E_1=2k_2I_n-2k_d\hat b-\frac{b^{\mathsf{T}}b}{2\underline b}.
\end{align*}
We next proceed to show that $D_1>0$,~ $E_1>0$ and
\begin{align*}
\lambda_{\min}(D_1)\lambda_{\min}(E_1)> \|B_1\|^2.
\end{align*}

First, from $\|a\|\le L_1$ and $\|b\|\le L_2$,
 we know
 $$\|\hat a\|=\left\|\left(a+a^{\mathsf{T}}\right)/2\right\|\le L_1,~~
 \|\hat b\|=\left\|\left(b+b^{\mathsf{T}}\right)/2\right\|\le L_2.$$ Besides, recall  $k_1=(k_p^2-2k_ik_d)\underline b,~ k_2=k_d^2\underline b-k_p$, it can be obtained that
\begin{align}
&\lambda_{\min} (D_1)\geq 2(k_p^2-2k_ik_d)\underline b-2k_pL_1-\frac{L_1^2}{2\underline b},\\
&\lambda_{\min} (E_1) \geq 2(k_d^2\underline b-k_p)-2k_dL_2-\frac{L_2^2}{2\underline b},\\
&\|B_1\|\le k_pL_2+k_dL_1+\frac{L_1L_2}{2\underline b}.
\end{align}
Since $(k_p,k_i,k_d)\in\Omega_{pid}$, it can be deduced that
\begin{align*}
k_p^2>&k_p^2-2k_ik_d>\bar k=(L_1+L_2)(k_p+k_d)/\underline b\\
>&(L_1+L_2)k_p/\underline b,
\end{align*}
which implies
$k_p>(L_1+L_2)/\underline b.$
Similarly, from $k_d^2>\bar k+k_p/\underline b>\bar k$, it can be verified that
$k_d>(L_1+L_2)/\underline b.$
Hence, the following inequalities hold:
\begin{align}\label{12}
\lambda_{\min} (D_1)\geq &2(k_p^2-2k_ik_d)\underline b-2k_pL_1-L_1^2/(2\underline b)\nonumber\\
\geq & 2(k_p+k_d)(L_1+L_2)-2k_pL_1-L_1^2/(2\underline b)\nonumber\\
= &2k_pL_2+2k_d(L_1+L_2)-L_1^2/(2\underline b)\nonumber\\
\geq &k_pL_2+k_d L_1+k_d(L_1+L_2)-L_1^2/(2\underline b)\nonumber\\
\geq& k_pL_2+k_d L_1+(L_1+L_2)^2/\underline b-L_1^2/(2\underline b)\nonumber\\
> & k_pL_2+k_dL_1+L_1L_2/(2\underline b)
\geq\|B_1\|.
\end{align}
Similarly, it can be deduced that
\begin{align}\label{16}
\lambda_{\min} (E_1)\geq & 2(k_d^2\underline b-k_p)-2k_dL_2-L_2^2/(2\underline b)\nonumber\\
\geq & 2(k_p+k_d)(L_1+L_2)-2k_dL_2-L_2^2/(2\underline b)\nonumber\\
> &k_pL_2+k_d L_1+k_p(L_1+L_2)-L_2^2/(2\underline b)\nonumber\\
\geq& k_pL_2+k_d L_1+(L_1+L_2)^2/\underline b-L_2^2/(2\underline b)\nonumber\\
\geq & k_pL_2+k_dL_1+L_1L_2/(2\underline b)
\geq \|B_1\|.
\end{align}
Combine (\ref{12}) with (\ref{16}), we know that $$\lambda_{\min}(D_1)\lambda_{\min}(E_1)> \|B_1\|^2.$$
By Lemma 6.2, we conclude that $\begin{bmatrix}D_1&B_1\\B_1^{\mathsf{T}}&E_1\end{bmatrix}>0$, i.e.  $E-B^{\mathsf{T}}D^{-1}B>0$, which implies $Q_0>0$.

\textbf{Step 5}. From Step 4, we know that $Q_0>0$ for all $\|a\|\le L_1$ and $\|b\|\le L_2$.
From the expression of $Q_0$ in (\ref{q00}), it is easy to see $Q_0$ is a continuous function of the variables $a$ and $b$. Notice that, the matrices $a$ and $b$ vary on a compact set, namely $\|a\|\le L_1$ and $\|b\|\le L_2$.
From Lemma 6.3, there exists $\alpha>0$(which depends on $(k_p,k_i,k_d,L_1,L_2,\underline b)$ only), such that
$Q_0\geq \alpha I_{3n}$ for all $\|a\|\le L_1$ and $\|b\|\le L_2$. By Step $3$, we know that $Q\geq Q_0\geq \alpha I_{3n}$, this implies that $A^{\mathsf{T}}P+PA\le -\alpha I_{3n}$ for all $a$, $b$ and $\theta $ with $\text{Sym}[\theta]\geq \underline bI_n$, $\|a\|\le L_1$ and $\|b\|\le L_2$.
 \hfill$\square$

\textbf{Proof of Proposition 4.3}
Since $(k_p,k_d)\in \Omega_{pd}$, then we know that $k_p>0$, $k_d>0$ and $k_d^2\underline b-k_p>\bar k>0$. Therefore, we have $2k_pk_d\underline b>0$, and
 $2k_pk_d^2\underline b-k_p^2=k_p(2k_d^2\underline b-k_p)>k_p(k_d^2\underline b-k_p)>0.$
Hence, the matrix $P$  defined by (\ref{P12}) is positive definite.
Next, we  show that $PA+A^{\mathsf{T}}P$ is negative definite.
Denote
$$Q\overset{\triangle}{=}-(PA+A^{\mathsf{T}}P)=\begin{bmatrix}Q_1&Q_{12}\\Q_{12}^{\mathsf{T}}&Q_{2}\end{bmatrix},$$ then by some calculations, we have
 $$Q_1=2k_p^2\text{Sym}[\theta]-k_p(a+a^{\mathsf{T}}),$$ $$~Q_{12}=2k_pk_d(\text{Sym}[\theta]-\underline bI_n)-k_pb-k_da^{\mathsf{T}}$$ and $$Q_2=2k_d^2\text{Sym}[\theta]-k_d(b+b^{\mathsf{T}})-2k_pI_n.$$

Notice that $\text{Sym}[\theta] \geq \underline b I_n$, similar  to Step 3 in the proof of Proposition 4.2, it is not difficult to obtain
\begin{align}\label{q2}Q\geq \begin{bmatrix}
2k_p^2\underline b I_n-k_p(a+a^{\mathsf{T}})&-k_pb-k_da^{\mathsf{T}}\\
-k_pb^{\mathsf{T}}-k_da&2(k_d^2\underline b-k_p) I_n-k_d(b+b^{\mathsf{T}})
\end{bmatrix}.\end{align}
For $z_1,z_2\in\mathbb{R}^n$, it follows from (\ref{q2}) that
\begin{align}\label{V}
&\left[z_1^{\mathsf{T}},z_2^{\mathsf{T}}\right]Q
\begin{bmatrix}
z_1\\z_2
\end{bmatrix}\nonumber\\\geq&
2k_p^2\underline bz_1^{\mathsf{T}} z_1-k_pz_1^{\mathsf{T}}(a+a^{\mathsf{T}})z_1- 2k_pz_1^{\mathsf{T}}bz_2\nonumber\\
&-2k_dz_1^{\mathsf{T}}a^{\mathsf{T}}z_2+2(k_d^2\underline b-k_p) z_2^{\mathsf{T}} z_2-k_dz_2^{\mathsf{T}}(b+b^{\mathsf{T}})z_2\nonumber\\
\geq& 2(k_p^2\underline b-k_pL_1)z_1^{\mathsf{T}}z_1-2(k_pL_2+k_dL_1)\|z_1\|\|z_2\|\nonumber\\
&+2(k_d^2\underline b-k_p-k_dL_2)z_2^{\mathsf{T}}z_2\nonumber\\
 \geq &  \left(2k_p^2\underline b-2k_pL_1-(k_pL_2+k_dL_1)\right)z_1^{\mathsf{T}} z_1\nonumber\\
 &+\left(2k_d^2\underline b-2k_p-2k_dL_2-(k_pL_2+k_dL_1)\right)z_2^{\mathsf{T}}z_2 \nonumber\\
\geq & \left(2k_p^2-2\bar k\right)\underline bz_1^{\mathsf{T}} z_1+\left(2k_d^2\underline b-2k_p-2\bar k\underline b\right)z_2^{\mathsf{T}}z_2,
\end{align}
where $\bar k=(k_p+k_d)(L_1+L_2)/\underline b$.

Since $(k_p,k_d)\in \Omega_{pd}$, we know that $k_p^2>\bar k$, $k_d^2-k_p/\underline b>\bar k,$ and therefore $Q>0.$
Notice that the RHS of (\ref{V}) does not depend on the matrices $a$, $b$ and $\theta$. Hence, if we denote  $$\beta=2\min\left\{\left(k_p^2-\bar k\right)\underline b,~k_d^2\underline b-k_p-\bar k\underline b\right\},$$
then $Q\geq \beta I_{2n}$(and thus $PA+A^{\mathsf{T}}P\le -\beta I_{2n}$), for all $n\times n$ matrices $a$, $b$ and $\theta$ with $\|a\|\le L_1$, $\|b\|\le L_2$ and $\text{Sym}[\theta]\geq \underline bI_n >0$.   \hfill$\square$

\textbf{Proof of Proposition 4.4} Since $(k_p,k_i)\in \Omega_{pi}$, then it is easy to obtain $k_p>0$, $k_i>0$ and $k_p^2\underline b-k_i>0$. Therefore, we have
 $2k_pk_i\underline b>0$ and $$ 2k_ik_p^2\underline b-k_i^2=k_i(2k_p^2\underline b-k_i)>k_i(k_p^2\underline b-k_i)>0.$$
Hence, the matrix $P$ defined by (\ref{P1}) is positive definite.

Next, we need to show that $PA+A^{\mathsf{T}}P$ is negative definite, i.e. $-(PA+A^{\mathsf{T}}P)>0$.
By some simple calculations, we have
\begin{align*}
Q\overset{\triangle}{=}-(PA+A^{\mathsf{T}}P)=\begin{bmatrix}Q_1&Q_{12}\\Q_{12}^{\mathsf{T}}&Q_{2}
\end{bmatrix},\end{align*}
where
\begin{align*}
&Q_1~~\!=~2k_i^2\text{Sym}[\theta],\\
&Q_{12}~=~2k_ik_p(\text{Sym}[\theta]-\underline bI_n)-k_ia,\\
&Q_2~~\!=~2k_p^2\text{Sym}[\theta]-k_p\left(a+a^{\mathsf{T}}\right)-2k_iI_n.
\end{align*}
First, notice that $\text{Sym}[\theta] \geq \underline b I_n$, similar to Step 3 in the proof of Proposition 4.2, it is not difficult to obtain
$$Q\geq \begin{bmatrix}
2k_i^2\underline b I_n&-k_ia\\
-k_ia^{\mathsf{T}}&~~~2k_p^2\underline b I_n-k_p(a+a^{\mathsf{T}})-2k_iI_n
\end{bmatrix}\overset{\triangle}{=}Q_0.$$
Next, we prove that  $Q_0>0$, for all $n\times n$ matrix $a$ with $\|a\|\le L$.

For $z_1,z_2\in\mathbb{R}^n$, it is easy to get
\begin{align}\label{V111}
&[z_1^{\mathsf{T}},~z_2^{\mathsf{T}}]Q_0\begin{bmatrix}z_1\\z_2\end{bmatrix}=
2k_i^2\underline b z_1^{\mathsf{T}} z_1+2k_p^2\underline bz_2^{\mathsf{T}}z_2\nonumber\\&
-k_pz_2^{\mathsf{T}}(a+a^{\mathsf{T}})z_2-2k_iz_2^{\mathsf{T}}z_2
- 2k_iz_1^{\mathsf{T}}az_2\nonumber\\
\geq&2k_i^2\underline b z_1^{\mathsf{T}} z_1+(2k_p^2\underline b-2k_p L-2k_i)z_2^{\mathsf{T}}z_2
- 2k_iLz_1^{\mathsf{T}}z_2.
\end{align}
From (\ref{V111}), it follows that
\begin{align}\label{V121}Q_0\geq Q_1\overset{\bigtriangleup}{=}\begin{bmatrix}2k_i^2\underline b &-k_i L \\
 -k_i L &~~2(k_p^2\underline b-k_p L-k_i)\end{bmatrix}\otimes I_n. \end{align}
Recall that $(k_p,k_i)\in\Omega_{pi}$, i.e., $4(k_p^2\underline b-k_pL-k_i)\underline b>L^2$,  it is easy to see that $Q_1$ is positive definite. Moreover, notice that $Q_1$ defined in (\ref{V121}) does not depend on the matrices $a$ and $\theta$. Hence, there exists $\gamma>0$($\gamma$ can be chosen as the smallest eigenvalue of $Q_1$), such that $Q_0\geq \gamma I_{2n}$(and thus $PA+A^{\mathsf{T}}P\le -Q_0\le -\gamma I_{2n}$), for all $n\times n$ matrices $a$ and $\theta$ with $\|a\|\le L$  and $\text{Sym}[\theta] \geq \underline bI_n >0$.
 \hfill$\square$

\subsection{Proof of Remark 3.1}\label{re1}
We first show that  $\Omega_{pid}$ defined by (\ref{pid}) is an open and unbounded subset in $\mathbb{R}^3$. The open property follows immediately from the definition (\ref{pid}).
Let $k_i>0$, $k_p=k_d=k\geq 2k_i+[2(L_1+L_2)+1]/\underline b$. Then, it can be obtained that
\begin{align*}
k_p^2-2k_ik_d-\bar k=&k^2-2k_ik-2(L_1+L_2)k/\underline b\\
=&k\left(k-2k_i-2(L_1+L_2)/\underline b\right)\\
\geq &k/\underline b>0,
\end{align*}
and
\begin{align*}
k_d^2-k_p/\underline b-\bar k=&k^2-k/\underline b-2(L_1+L_2)k/\underline b\\
=&k\left(k-1/\underline b-2(L_1+L_2)/\underline b\right)\\
\geq &2k_ik>0,
\end{align*}
which implies that the following set is a subset of  $\Omega_{pid}$:
$$\Big\{(k_p,k_i,k_d)\Big| k_i>0,~ k_p=k_d\geq 2k_i+[2(L_1+L_2)+1]/\underline b\Big\}.$$  Thus,  $\Omega_{pid}$ is unbounded and not empty.

Next, we show that $\Omega_{pid}$ is a semi-cone in the following  sense:
$(k_p,k_i,k_d)\in\Omega_{pid}\Rightarrow\alpha(k_p,k_i,k_d)\in\Omega_{pid},~~\forall \alpha\geq 1.$
Suppose that  $(k_p,k_i,k_d)\in\Omega_{pid}$ and $\alpha\geq 1$, we need to show   $(k_p', k_i', k_d')=\alpha(k_p,k_i,k_d)\in\Omega_{pid}$. Let us define $\bar k'=(L_1+L_2)(k_p'+k_d')/\underline b$. Then it is easy to obtain
\begin{align}\label{ty}
k_p'^2-2k_i'k_d'-\bar k'=&\alpha^2(k_p^2-2k_ik_d)-\alpha \bar k\nonumber\\
\geq & \alpha^2\left(k_p^2-2k_ik_d- \bar k\right)>0.
\end{align}
Therefore, the first inequality in the definition of $\Omega_{pid}$ is satisfied.
Next, it is easy to obtain the  following inequalities:
\begin{align}\label{tz}
k_d'^2-k_p'/\underline b-\bar k'=&\alpha^2k_d^2-\alpha k_p/\underline b-\alpha \bar k\nonumber\\
\geq &\alpha^2\left(k_d^2- k_p/\underline b- \bar k\right)>0.
\end{align}
Thus,  the second inequality is also verified. Therefore, $\alpha(k_p,k_i,k_d)\in\Omega_{pid}$ for all $\alpha\geq 1$. \hfill$\square$


\begin{thebibliography}{20}
	\expandafter\ifx\csname natexlab\endcsname\relax\def\natexlab#1{#1}\fi
	\expandafter\ifx\csname url\endcsname\relax
	\def\url#1{\texttt{#1}}\fi
	\expandafter\ifx\csname urlprefix\endcsname\relax\def\urlprefix{URL }\fi
	
\bibitem[{\AA}str{\"o}m \& H{\"a}gglund, 2000]{As}
    {\AA}str{\"o}m, K. J., \& H{\"a}gglund, T. (2000). The future of PID control. {\em Control Engineering Practice}, 9(11): 1163-1175.
\bibitem[{\AA}str{\"o}m \&  H{\"a}gglund, 1995]{Astrom19955}
{\AA}str{\"o}m, K. J., \&  H{\"a}gglund, T. (1995).
\newblock {\protect{PID Controllers: Theory, Design and Tuning}}.
\newblock  Research Triangle Park, NC: Instrument society of America.
\bibitem[Samad, 2017]{samad} Samad, T. (2017).  A Survey on Industry Impact and Challenges Thereof. {\em IEEE Control Systems Magazine}, 37(1): 17-18.
\bibitem[Ho \& Lin, 2003]{Ho2003} Ho, M., \& Lin,  C. (2003). PID controller design for robust performance. {\em IEEE Transactions on Automatic Control},  48(8):1404--1409.
\bibitem[Keel \&   Bhattacharyya, 2008]{Keel2008}
Keel, L. H., \&   Bhattacharyya, S. P. (2008).
Controller synthesis free of analytical models: Three term
  controllers. {\em IEEE Transactions on Automatic Control}, 53(6):1353--1369.

\bibitem[Silva {\em et al.},  2005]{Silva}
Silva, G. J.,  Datta, A., \&  Bhattacharyya, S. P. (2005). PID controllers for time-delay systems. Boston: Birkhäuser.

\bibitem[Ziegler \& Nichols, 1993]{ziegler1942}
Ziegler, J. G., \& Nichols, N. B. (1993).
\newblock{\protect{Optimum Settings for Automatic Controllers.}}
 {\em Journal of Dynamic Systems Measurement and Control-transactions of The Asme}, 220-222.

\bibitem[Odwyer, 2006]{o2006pi}
Odwyer, A. (2006). PI and PID controller tuning rules: an overview and personal perspective. in Proc. IET Irish Signals Syst. Conf., Dublin, Ireland, 2006: 161-166.

\bibitem[Takegaki \& Arimoto, 1981]{Takegaki} Takegaki, M., \&  Arimoto, S. (1981). A new feedback method for dynamic control of manipulators. {\em Journal of Dynamic Systems Measurement and Control-transactions of The Asme}, 103(2), 119-125.

\bibitem[Alvarez-Ramirez {\em et al.}, 2002]{Alvarez}
Alvarez-Ramirez, J., Kelly, R., \&  Cervantes, I. (2002). Semiglobal stability of saturated linear PID control for robot manipulators, {\em Automatica}, 39 (6), 989–995.

\bibitem[Besancon, 2000]{Besancon} Besancon, G. (2000). Global output feedback tracking control for a class of Lagrangian systems. {\em Automatica}, 36 (12), 1915-1921.

\bibitem[Su {\em et al.}, 2010]{Su}
Su, Y., Muller, M. C.,  Zheng, C. (2010).  Global asymptotic saturated PID control for robot manipulators.
{\em IEEE Transactions on Control Systems Technology},  18 (6), 1280–1288.

\bibitem[Borja {\em et al.}, 2021]{Borja}
Borja, P.,  Ortega, R., \&  Scherpen, J. M. A. (2021). New Results on Stabilization of Port-Hamiltonian Systems via PID Passivity-Based Control. {\em IEEE Transactions on Automatic Control}, vol. 66, no. 2, pp. 625-636,  doi: 10.1109/TAC.2020.2986731.

\bibitem[Killingsworth \& Krstic, 2006]{Killingsworth2006}
 Killingsworth,  N. J., Krstic, M. (2006).
 \protect{PID tuning using extremum seeking: online, model-free performance optimization, {\em IEEE Control Systems Magazine},  26(1): 70-79.}

\bibitem[Khalil, 2000]{Khalil} Khalil, H. K. (2000). Universal integral controllers for minimum-phase nonlinear systems. {\em IEEE Transactions on Automatic Control},  45(3):490-494.

\bibitem[Chang {\em et al.}, 2002] {Chang} Chang, W. D.,  Hwang, R. C.,  Hsieh, J. G. (2002). A self-tuning PID control for a class of nonlinear systems based on the Lyapunov approach. {\em Journal of Process Control}, 12(2):233-242.

\bibitem[Romero {\em et al.}, 2018]{Romero} Romero, J. G.,  Donaire, A., Ortega, R., Borja, P. (2018). Global stabilisation of underactuated mechanical systems via PID passivity-based control. {\em Automatica}, 96(96): 178-185.
\bibitem[Kelly, 1998]{kelly}Kelly, R. (1998).  Global positioning of robot manipulators via PD control plus a class of nonlinear integral actions, {\em IEEE Transactions on Automatic Control}, vol. 43, no. 7, pp. 934-938.
\bibitem[Brogliato {\em et al.}, 2007]{Brogliato} Brogliato, B., Lozano, R., Maschke, B.,  Egeland, O. (2007). Dissipative systems analysis and control. Theory and Applications, 2.
\bibitem[Guo, 2020]{lguo}
 Guo, L. (2000). Feedback and Uncertainty: Some Basic Problems and Results, {\em Annual Reviews in Control}, https://doi.org/10.1016/j.arcontrol.2020.04.001.
 \bibitem[Jiang {\em et al.}, 2015]{jiang2015} Jiang, T.T.,  Huang, C.D., Guo, L. (2015). Control of uncertain nonlinear systems based on both observers and estimators. {\em Automatica}, vol.59, pp.35-47.

\bibitem[Guo \& Zhao, 2021]{guo2021} Guo, L., \& Zhao, C. (2021). Control of nonlinear uncertain systems by extended PID with differential trackers. {\em Communications in Information and Systems}, vol. 21, No.3, 415-440.

\bibitem[Krstic, 2017]{krstic2017}
 Krstic, M. (2017). On the applicability of PID control to nonlinear second-order systems. {\em National Science Review}, vol 4, issue 5, page 668.
 \bibitem[Xie \& Guo, 2000]{xie2000} Xie, L. L., \& Guo, L. (2000). How much uncertainty can be dealt with by feedback?.  {\em IEEE Transactions on Automatic Control},  45(12): 2203-2217.

\bibitem[Xue \& Guo, 2002]{xue2002} Xue, F., \& Guo, L. (2002). On limitations of the sampled-data feedback for nonparametric dynamical systems, J. Systems Science and Complexity, 15(3):225-250.

\bibitem[Zhao \& Guo, 2017]{Zhao2017} Zhao,  C., \& Guo, L. (2017). PID controller design for second order nonlinear uncertain systems. {\em Science China Information Sciences}, 60(2): 022201,  doi: 10.1007/s11432-016-0879-3.
\bibitem[Zhang \&  Guo, 2019]{Zhang2019}Zhang, J. K., \&  Guo, L.  (2019). Theory and design of PID controller for nonlinear uncertain systems. {\em IEEE Control Systems Letters}, 3(3):643--648.
\bibitem[Zhao \& Guo, 2021]{zhao2020}  Zhao, C., \& Guo, L. (2021).  Control of nonlinear uncertain systems by extended PID, {\em IEEE Transactions on Automatic Control}, vol.66, no.8, 3840-3847.
\bibitem[Zhao  \&  Guo, 2018]{Zhao2018} Zhao, C., \&  Guo, L. (2018). PID control for a class of non-affine uncertain systems. 37th Chinese Control Conference,  doi:10.23919/chicc.2018.8483587.
 \bibitem[Zhao \& Guo, 2017b]{zhao2016}Zhao, C., \& Guo, L. (2017b). On the capability of PID control for nonlinear uncertain systems. IFAC-PapersOnLine,  50(1): 1521-1526.
\bibitem[Boskovic {\em et al.}, 2004]{bos} Boskovic, J. D., Chen, L., Mehra, R. K. (2004).   Adaptive control design for nonaffine models arising in flight control. {\em Journal of Guidance Control and Dynamics},  27(2): 209-217.

\bibitem[Zhong {\em et al.}, 2021]{Zhong} Zhong, S., Huang, Y.,  Guo, L. (2021). An ADRC-based PID tuning rule, International Journal of Robust and Nonlinear Control, 2021; doi: 10.1002/rnc.5845.

\bibitem[Huang {\em et al.}, 2001]{huang2001}
    Huang, Y., Xu, K., Han, J., \& Lam, J. (2001). Flight control design using extended state observer and non-smooth feedback. In Proceedings of the 40th IEEE Conference on Decision and Control(Cat. No. 01CH37228), vol. 1, pp. 223-228).
\bibitem[Gordon, 1972]{gordon1972}
Gordon, W. B. (1972) On the diffeomorphisms of Euclidean space. The American Mathematical Monthly, 1972, 79(7): 755-759.

\bibitem[Ruzhansky \& Sugimoto, 2015]{Ruzhansky2015}
Ruzhansky, M., Sugimoto, M. (2015) On global inversion of homogeneous maps. Bulletin of Mathematical Sciences, 2015, 5(1): 13-18.
\bibitem[Fe{\ss}ler, 1995]{fessler1995proof}
Fe{\ss}ler, R. (1995). A proof of the two-dimensional Markus-Yamabe stability conjecture and
  a generalization.
{\em Ann Polon Math}, 62:45--74.
\bibitem[Horn \&  Johnson, 1985]{horn}Horn, R. A.,  \& Johnson, C. R. (1985). Matrix analysis. Cambridge University Press, 1985.
\end{thebibliography}
\end{document}